\documentclass[12pt]{article}
\usepackage{amsfonts,amssymb,version,amsmath,graphicx,multicol}
\allowdisplaybreaks

\def\mm{\mathfrak m}

\def\P{{\cal P}}

\def\R{\mathbb R}

\def\T{{\mathbb T}}

\def\x{\dot{x}}

\mathchardef\mhyphen="2D

\mathchardef\mhyphen="2D

\let\hra\hookrightarrow

\let\un\underline
\let\upa\uparrow

\let\pa\partial

\newtheorem{theorem}{Theorem}[section]

\newtheorem{lemma}[theorem]{Lemma}

\def\stm{\refstepcounter{theorem}\paragraph{\thetheorem}}

\makeatletter \def\l@section{\@dottedtocline{1}{0em}{1.2em}} \makeatother

\begin{document}

\centerline{\Large\bf Curvature, torsion and the quadrilateral gaps}

\medskip


\bigskip

\centerline{\bf Nitin Nitsure}


\begin{abstract} 
For a manifold with an affine connection, we prove formulas which 
infinitesimally quantify the gap in a certain naturally defined open geodesic 
quadrilateral associated to a pair of tangent vectors $u$, $v$ at a 
point of the manifold. 
We show that the 1st order infinitesimal obstruction to the 
quadrilateral to close is always zero, the 2nd order infinitesimal 
obstruction to the quadrilateral to close is $-T(u,v)$ where $T$ is 
the torsion tensor 
of the connection, and if $T = 0$ then the 3rd 
order infinitesimal obstruction to the quadrilateral to close is 
$(1/2)R(u,v)(u+v)$ in terms of the curvature tensor of the connection. 
Consequently, the torsion of the connection, and if 
the torsion is identically zero then also the curvature of the connection, 
can be recovered uniquely from knowing all the quadrilateral gaps. 
In particular, this answers a question of Rajaram Nityananda about the 
quadrilateral gaps on a curved Riemannian surface. The angles
of $3\pi/4$ and $-\pi/4$ radians feature prominently in the answer, 
along with the value of the Gaussian curvature. This article is 
essentially self-contained, and written in an expository style.

\end{abstract}

\pagestyle{empty}


%

\section{Introduction}

Let $M$ be a surface with a Riemannian metric, and let $P_0$ be a point on $M$.
Let $u\in T_{P_0}M$ be a unit tangent vector. 
Now consider the following journey on $M$. 
To begin with, choose the geodesic starting at $P_0$ in the direction
given by $u$, and travel along it 
for a distance $s$ to arrive at a point $P_1(s)$.  
Next, turn left in $\pi/2$ radians, and travel along the geodesic 
in that direction for a distance $s$, to arrive at a point $P_2(s)$. 
Repeat this twice more: turn left at $P_2(s)$ in $\pi/2$ radians 
and travel along the geodesic in that direction for a distance $s$ to arrive
at a point $P_3(s)$, and then turn left at $P_3(s)$ in $\pi/2$ radians and 
travel along the geodesic 
in that direction for a distance $s$ to finally arrive at a point $P_4(s)$. 
This defines an open rectangle on $M$, with vertices $P_0,\ldots, P_4$,
whose legs are the above four geodesic segments.
We assume that we have chosen an orientation for $M$ around $P_0$, so that 
the left turns make unambiguous sense.

If the surface $M$ is flat, then for small values of $s$, we will be back
at the starting point, that is, our journey will be along a geodesic 
quadrilateral with $P_4(s) = P_0$. But if $M$ is not flat, then we do not come
back, that is, $P_4(s) \ne P_0$. One may say that the quadrilateral 
does not close, as it has a gap. Rajaram Nityananda asked for a 
quantification of this gap for a small $s$.

We have another geodesic starting at $P_0$ with initial direction
$v\in T_{P_0}M$ where $(u,v)$ is a right-handed orthonormal basis for 
$T_{P_0}M$. Travel for a distance $s$ along it, to come to a point 
$Q_1(s)$. Turn right at $Q_1(s)$ in $\pi/2$ radians and 
travel along the geodesic 
in that direction for a distance $s$ to arrive at a point $Q_2(s)$.
Again, if $M$ is flat, then $Q_2(s) = P_2(s)$ for a small $s$.
Otherwise, we have an open rectangle with successive vertices 
$Q_2, Q_1,P_0,P_1,P_2$ joined by geodesic segments of length $s$.
Again, one can ask what is the gap between $P_2(s)$ and $Q_2(s)$.


\pagestyle{myheadings}
\markright{Nitin Nitsure: 
Curvature, torsion and the quadrilateral gaps}


Here is our answer. We show that
$$\lim_{s\to 0} {P_4(s) - P_0 \over s^3} = 
\lim_{s\to 0} {Q_2(s) - P_2(s) \over s^3} = 
{\kappa(P_0) \over 2}
(u-v)$$
where $\kappa(P_0)$ denotes the Gaussian curvature of $M$ at $P_0$.

{
\begin{center}
\includegraphics[scale=.08]{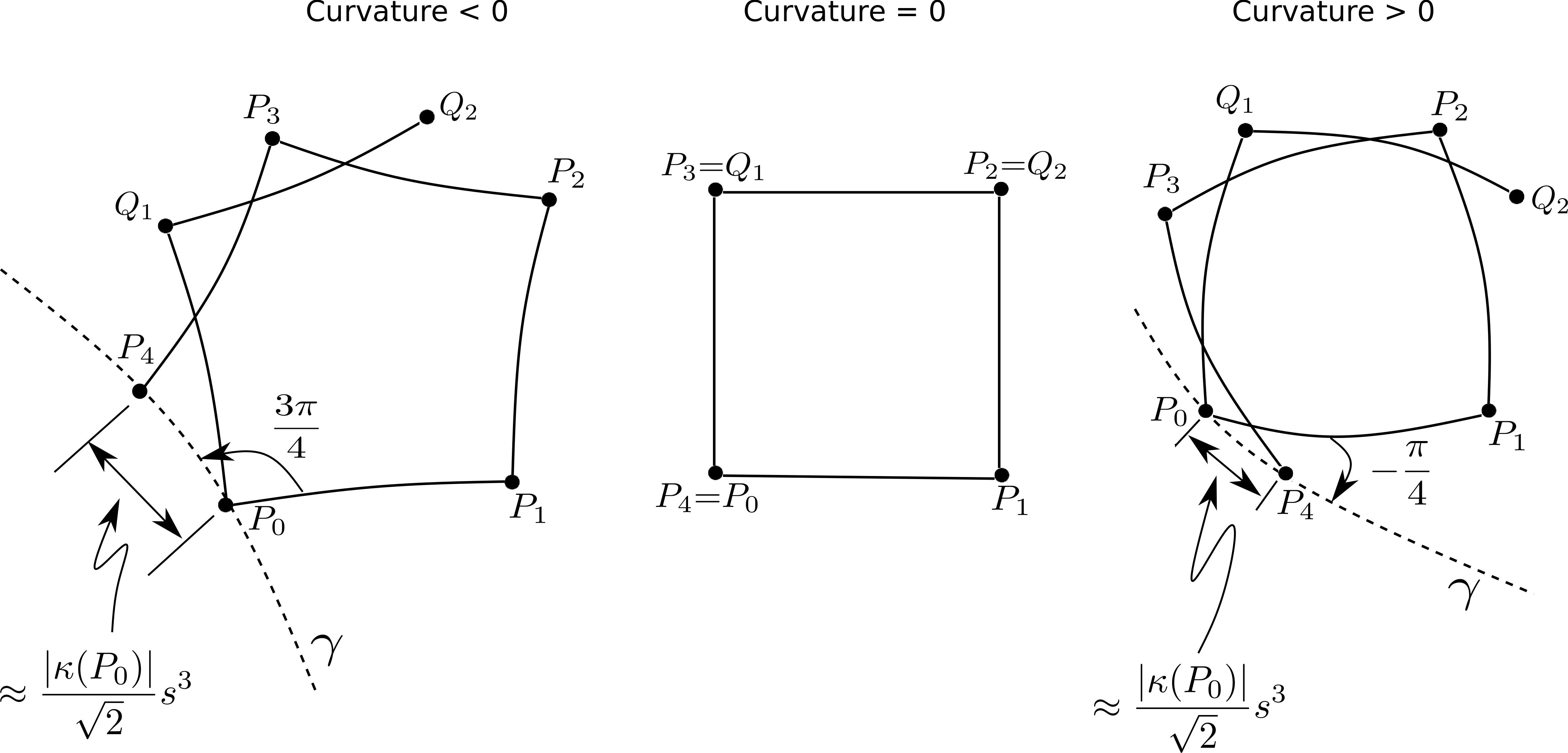}
\end{center}
}

The difference $P_4(s) - P_0$ or $Q_2(s) - P_2(s)$ 
between any pair of points of $M$ 
can be understood in terms of a smooth 
embedding of a neighbourhood of $P_0$ into an affine space $\R^n$,
which makes the difference a vector in the vector space $\R^n$. 
In particular, the points $P_4(s)$ trace a curve $\gamma(\lambda)$, 
parameterized by $\lambda = s^3$, and this has the tangent vector 
$(\kappa(P_0)/2)(u-v)$ at $\lambda =0$.

Equivalently, 
for any smooth function $f$ defined in a neighbourhood of $P_0$, we have
$$\lim_{s\to 0} {f(P_4(s)) - f(P_0) \over s^3} = 
\lim_{s\to 0} {f(Q_2(s)) - f(P_2(s))  \over s^3} = 
{\kappa(P_0) \over 2}
(u-v)(f).$$
In words, the sizes of both the gaps, that is, the distances
from $P_0$ to $P_4(s)$ or from $P_2(s)$ to $Q_2(s)$ 
regarded as functions of $s$, 
are equal to ${1\over \sqrt{2}}|\kappa(P_0)|s^3$ upto 3rd order terms in $s$
(that is, modulo $s^4$). 
Secondly, upto 3rd order terms in $s$,
the gaps make an angle of $-\pi/4$ radians with the first leg of the 
quadrilateral when $\kappa(P_0) > 0$, and 
the opposite angle of $3\pi/4$ radians when $\kappa(P_0) < 0$, 
regardless of the magnitude $|\kappa(P_0)|$ of the curvature.

The above result may be compared with 
the formula of Bertrand and Puiseux, which gives the 
circumference $C(r)$ of an infinitesimally small circle on $M$ 
centered at $P_0\in M$ with geodesic radius $r$. 
This formula says that 
the deviation of $C(r)$ from the Euclidean value 
$2\pi r$ is infinitesimally small of 3rd order in $r$. More precisely, 
$$\lim_{r\to 0} {2\pi r - C(r) \over r^3} = {\pi\kappa(P_0)\over 3},$$ 
that is,  
$2\pi r - C(r) = {1\over 3}\pi\kappa(P_0) r^3$
upto 3rd order terms in $r$. 

We directly verify in Section 3 the result for the sphere $x^2+y^2+z^2 = r^2$
in the Euclidean space $\R^3$ (which has constant positive curvature
$1/r^2$) and 
for the hyperboloid $t^2 - x^2 - y^2 = r^2$ in the Minkowski
space $\R^{1,2}$ (which has constant negative curvature
$-1/r^2$) 
by respectively using orthogonal or Lorentz transformations.
We sketch a heuristic argument in Section 3, 
based on Riemann normal coordinates, 
to go from the case of constant curvature to that of arbitrary metrics
on surfaces. We will not make this argument rigorous, as 
the formula
$\lim_{s\to 0}(P_4(s) - P_0)/s^3 = (\kappa(P_0)/2)(u-v)$ is just the 
$2$-dimensional case of the more general Theorem \ref{geodesic gap},
which we now formulate.

Let $M$ be a differential manifold of any dimension $d$, 
equipped with an affine connection $\nabla$, 
that is, a connection on the tangent bundle $TM$.
As a special case,
$M$ may be a (pseudo-)Riemannian manifold, and $\nabla$ the 
connection induced by the given metric (which is the unique
symmetric connection for which the covariant derivative of the 
metric is identically zero). 
We now generalize Nityananda's question to this general setup,  
where notions such as length and perpendicularity are not available. 
To any ordered triple $(P,u,v)$ consisting of a point $P\in M$ and a pair
of tangent vectors $u,v \in T_PM$, and any real number $s$ such that 
$|s|$ is small enough (depending on $(P,u,v)$), we associate
a new such triple $\T_s(P,u,v) = (P',u',v')$ constructed as follows.
Let $\gamma$ be the unique affinely parameterized geodesic starting at $P$
with initial tangent vector $u$, so that $\gamma(0) = P$ and 
$\dot\gamma(0) = u$.  
We define
$$P' = \gamma(s),~ u' = v(s) \mbox{ and } v' = - u(s)$$
where $u(s),v(s) \in T_{P'}M$ are the parallel transports 
of $u,v$ along $\gamma$. 
We will sometimes write the 
triple $(P',u',v')$ as $(P'(s), u'(s),v'(s))$
to make the dependence on $s$ explicit. As an example,
$\T_0(P,u,v) = (P'(0), u'(0),v'(0)) = (P,v,-u)$. 
The point $P$ will be called as the {\bf location} of the triple $(P,u,v)$.
We now begin with a point $P_0\in M$ and $u,v\in T_{P_0}M$ and apply
the operator $\T_s$ iteratively to define new points $P_1, P_2, \ldots$, 
where $P_1$ is the location of $\T_s(P,u,v)$, $P_2$ is the location of 
$\T_s(\T_s(P,u,v))$, and in general  
$P_n$ is the location of the triple $(\T_s)^n(P,u,v)$. 
Note that the operator $\T_s$ is invertible (but its two sided inverse 
$\T^{-1}_s$ is not equal to $\T_{-s}$ in general). We 
apply the operator $\T^{-1}_s$ iteratively to define new points 
$Q_1, Q_2, \ldots$, 
where $Q_1$ is the location of $T^{-1}_s(P,u,v)$, $Q_2$ is the location of 
$\T^{-1}_s(\T^{-1}_s(P,u,v))$, and so on. 
In these terms, the quadrilateral gap of Nityananda is the gap between 
$P_0$ and $P_4(s)$ (or between $P_2(s)$ and $Q_2(s)$) when $M$ is Riemannian 
of dimension $2$ and $(u,v)$ is an 
orthonormal basis for $T_{P_0}M$.
The Theorem \ref{geodesic gap} describes these gap infinitesimally.

\theorem\label{geodesic gap}
{\it Let $M$ be a smooth manifold equipped with an affine connection
$\nabla$. Then with notation as above, we have the following.

(1) If $P_0 \in M$ and  $u,v \in T_{P_0}M$,
then  
$$\lim_{s\to 0} {P_4(s) - P_0 \over s^2} = 
\lim_{s\to 0} {P_2(s) - Q_2(s)\over s^2} = -T(u,v) \in T_{P_0}M$$
where $T(u,v) = \nabla_uv - \nabla_vu - [u,v]$
is the torsion tensor of $\nabla$. 

(2) If $T\equiv 0$, that is, $\nabla$ is symmetric,
then 
$$\lim_{s\to 0} {P_4(s) - P_0\over s^3} = 
- \lim_{s\to 0} {P_2(s) - Q_2(s)\over s^3} 
= {1\over 2}R(u,v)(u+v) \in T_{P_0}M$$
where $R(u,v) = \nabla_u\nabla_v - \nabla_v\nabla_u - \nabla_{[u,v]}
\in End(T_{P_0}M)$ is the curvature tensor of $\nabla$. 
}

\rm

Our answer to Nityananda's question now
follows by taking $(u,v)$ to be 
an orthonormal basis
for $T_{P_0}M$ when $M$ is a surface with a Riemannian metric and 
$\nabla$ is the Riemannian connection, and noting that 
in this case, $T(u,v) =0$ and 
$R(u,v)$ is the skew-symmetric operator
$\left(
\begin{array}{cc}
0 & \kappa_0 \\
-\kappa_0 & 0
\end{array}
\right)$ 
in terms of the basis $(u,v)$, 
where $\kappa_0$ is the Gaussian curvature of $M$ at $P_0$.

Thus, the torsion tensor $T$ can be read off from the limit as $s\to 0$
of gaps$/s^2$. If the torsion is zero then
the elementary Lemma \ref{curvature identities} shows that 
the curvature tensor $R$ at a point $P_0$ is uniquely determined by
the function $(u,v)\mapsto R(u, v)(u+v)$, and so 
the curvature tensor $R$ can be uniquely recovered from the limit as $s\to 0$
of gaps$/s^3$. 

{\bf A heuristic argument.} Without knowing the Theorem \ref{geodesic gap},
one could have argued as follows. We can expand the 
gap $Q_2(s) - P_2(s)$ as a power series in $s$, to write 
$$Q_2(s) - P_2(s) = A_0(u,v) 
+ s A_1(u,v) + s^2 A_2(u,v) + s^3A_3(u,v)  \mbox{ mod }s^4.$$
The vectors $A_i(u,v)$ will have to be homogeneous polynomials 
of total degree $i$
in the variables $u,v$, as $A_i(u,v)$ is the coefficient of $s^i$,
and moreover, they should change sign when $u$ and $v$ are interchanged.
As we know that $Q_2(0) = P_2(0) = P_0$, we must have $A_0 =0$. 
The $A_i$'s are to be made from $\nabla$ alone. 
The obvious candidate for $A_1$ is $c(u-v)$ for some constant $c$.
But for the Euclidean space itself, the gap is zero, so we must have
$c=0$. This means $A_1(u,v) = 0$. An obvious candidate 
for $A_2(u,v)$ is $T(u,v)$, the torsion. A degree $3$ homogeneous candidate
for $A_3(u,v)$ is $R(u,v)u$, but that is not skew-symmetric. 
Applying the idempotent projector for skew-symmetry in the group ring 
$\R[S_2]$ of the permutation group $S_2$ to $R(u,v)u$, we get
the candidate ${1\over 2}(R(u,v)u - R(v,u)v)$ for $A_3$, which 
is homogeneous of degree $3$ and skew-symmetric in $u,v$. But this
equals ${1\over 2}R(u,v)(u+v)$, as $R(u,v)$ is skew-symmetric in $u,v$.
The above argument (made in hindsight) 
is only suggestive: for example, there could be constant numerical coefficients
$c_i$ that multiply the candidates $A_i$, though the examples in Section 2
would show us that there is no further coefficient that multiplies
${1\over 2}R(u,v)(u+v)$. 
The actual proof of Theorem \ref{geodesic gap} is given in 
Section 5, and it does not refer to the above argument.

In our discussion so far, we had chosen a pair of vectors $u,v\in T_{P_0}M$,
and we had parallel translated this pair. More generally, we can choose
a basis $(u_1,\ldots,u_d)$ for $T_{P_0}M$ where $d = \dim(M)$, and parallel 
transport it. This suggests that we should consider paths in the 
frame bundle $\pi: E \to M$ of $M$,
which is a principal $GL(d)$-bundle. 
Recall that a connection $\nabla$ on $TM$ naturally gives rise to a 
vector field $\xi$ on the 
total space of $TM$, whose flow is called as the geodesic flow of $\nabla$.
Similarly, $\nabla$ naturally gives rise to $d$ different horizontal 
vector fields
$\xi_1,\ldots, \xi_d$ on the total space of $E$. 
Here, recall that giving a connection $\nabla$ on $TM$ is equivalent to  
giving a rank
$d$ vector subbundle $D$ of $TE$, which is complementary to the kernel 
of the derivative map $\pi_* : TE \to \pi^* TM$ (the `vertical' subbundle)
at each point of $E$. For $y\in E$, an element of $T_yE$ is
called a horizontal vector if it lies in the fiber $D_y \subset T_yE$.  
If $x\in M$ and if $(u_1,\ldots, u_d)$ is a basis for $T_x M$, so that 
we can regard $y = (u_1,\ldots, u_d)$ as a point of $E$ lying over $x$,
then these horizontal vector fields $\xi_i$ are uniquely
defined by the requirement that 
$$\pi_*(y)\xi_i = u_i$$ 
where $\pi_*(y) : T_yE \to T_xM$ is the derivative of $\pi$. 
These vector fields have been considered in the book of K. Nomizu 
(see Chapter 3, section 1 of [1]) under the name `basic vector fields'.
With this definition, if we start at a point $P_0$ and parallel transport
a basis $y = (u_1,\ldots,u_d)$ along the affinely parameterized 
geodesic $\gamma$ through 
$P_0$ with $\gamma(0) = P_0$ and initial tangent vector $\dot\gamma(0) = u_i$
to the point $P_1= \gamma(s)$, and if $y_1 = (u_1(s), \ldots, u_d(s))$ 
is the basis of $T_{P_1}M$ obtained by parallel translating the basis 
$y= (u_1,\ldots,u_d)$ along $\gamma$, then $y_1\in E$ is exactly the point
obtained by following the flow on $E$ of the vector field $\xi_i$ for
parameter value $s$. This translates the question of geodesic 
quadrilateral gaps on $M$ into a question about non-commutativity of
flows of vector fields on the manifold $E$. The geodesic gap 
corresponding to $(P_0,u_m,u_n)$ will just be the projection
under $\pi : E\to M$ of the gap in $E$ corresponding
to the flows of vector fields $\xi_m$ and $\xi_n$. 
 
Let $\xi$ and $\eta$ be 
vector fields on a manifold $M$ and let $\varphi_s$ and $\psi_s$ 
be the corresponding flows. Given $P_0\in M$, for a small enough $s$ 
we can define the points 
$$P_1 = \varphi_s(P_0), ~ P_2 = \psi_s\circ\varphi_s(P_0), 
~ P_3 = \varphi_{-s}\circ\psi_s\circ\varphi_s(P_0), 
P_4 = \psi_{-s}\circ\varphi_{-s}\circ\psi_s\circ\varphi_s(P_0),$$
and 
$$
Q_1 = \psi_s(P_0) 
\mbox{ and }~ Q_2 = \varphi_s\circ\psi_s(P_0)
$$
in $M$. The points $P_0,P_1,P_2,P_3,P_4$ and $Q_{-2}, Q_{-1}, P_0,P_1,P_2$
are the successive vertices of open quadrilaterals, whose sides
are the flow lines of $\pm\xi$ and $\pm\eta$.  
The following is a well know fundamental result, that
interprets the bracket of two vector fields in terms of 
the gaps $P_4(0)-P_0$ and $P_2(s) - Q_2(s)$.

\lemma\label{meaning of bracket}
{\it With notation as above
$$\lim_{s\to 0} {P_4(s) - P_0\over s^2} = 
\lim_{s\to 0} {P_2(s) - Q_2(s)\over s^2} = [\xi,\eta]_{P_0} \in T_{P_0}M.$$
}
\rm

For the sake of completeness, we give a proof of the above in Section 3.
This proof introduces a certain trick which is important in the proof of 
Theorem \ref{geodesic gap}.

We now apply the above lemma to the flows of $\xi_m$ and $\xi_n$ on $E$, 
starting at a point $y=(u_1,\ldots,u_d)\in E$ over $x = P_0 \in M$.
The projection under $\pi: E\to M$ of the corresponding flow quadrilateral gaps 
in $E$ are the geodesic quadrilateral gaps in $M$
for the triple $(P_0,u_m,u_n)$. Thus, we need to compute the bracket 
$[\xi_m,\xi_n]$ in $E$ and its projection under $\pi_*: T_yE \to T_xM$. 
This is given by the following result
proved in Section 4, which follows from the definitions by an easy computation.
We expect this result -- or some equivalent form --
to be well known to experts. 

\theorem\label{frame flow bracket} {\it 
Let $\pi : E \to M$ be the frame bundle associated to $TM$.
Let $\xi_1, \ldots,  \xi_d$ be the frame flow fields on
the total space of $E$, associated to the affine 
connection $\nabla$ on $M$. Let $x \in M$ and let
$y = (u_1, \ldots, u_d)\in E$ with $\pi(y) = x$,
so that $\pi_*(y)(\xi_i) = u_i$ for all $i$. Then the following holds:

(1) The torsion of $\nabla$ is given by 
$T(u_i,u_j) = - \pi_*(y) [\xi_i,\xi_j]$. 

(2) In particular if  $\nabla$ is symmetric 
then $\pi_*[\xi_i,\xi_j] =0$, so
the bracket $[\xi_i,\xi_j]$ is a vertical vector field on $E$, 
hence can be naturally regarded as a section of the pullback 
to $E$ of 
the bundle $Ad(E) = \un{End}(TM)$ on $M$. 

(3) If $\nabla$ is symmetric, then the curvature of $\nabla$ is 
given as follows.  
For any $x\in M$ and $y = (u_1,\ldots, u_d)\in E$
with $\pi(y) =x$, we have   
$R(u_i,u_j)_x = -[\xi_i,\xi_j]_y$
as elements of $End(T_xM)$. 
As global sections of $\pi^*\un{End}(TM)$ over $E$,
we have the equality
$R(\pi_*\xi_i,\pi_*\xi_j) = -[\xi_i,\xi_j]$.
}
\rm

The fact that the geodesic gap for $(P_0,u_m,u_n)$
in $M$ is the projection under $\pi: E\to M$ of the flow gap
for $\xi_m,\xi_n$ in $E$, together with the Lemma \ref{meaning of bracket},
shows that the statement (1) of the 
Theorem \ref{frame flow bracket}
gives a proof of the statement (1) of the Theorem \ref{geodesic gap}.

The verticality of $[\xi_i, \xi_j]$ for symmetric connections 
shows that upto 2nd infinitesimal order in $s$,  
the geodesic quadrilateral gap on $M$ is zero for symmetric connections. 
In fact, we proved the Theorem \ref{frame flow bracket} before proving the 
Theorem \ref{geodesic gap}, and this
told us that the geodesic quadrilateral gap on $M$ for a symmetric connection 
is a phenomenon of 3rd or higher infinitesimal order in $s$. 
This is what motivated the 3rd order Taylor series calculation which is 
at the heart of the proof of the Theorem \ref{geodesic gap}.

In conclusion, we can therefore say that the 1st order 
infinitesimal obstruction to the quadrilaterals to close is always zero,
the 2nd order infinitesimal  
obstruction to the quadrilateral associated to $u,v\in T_{P_0}M$ to 
close is $T(u,v)$, and $\nabla$ is symmetric, then the 3rd 
order infinitesimal obstruction to this quadrilateral to close is 
$(1/2)R(u,v)(u+v)$. If $T \equiv 0$, then the 
Lemma \ref{curvature identities} shows
how to uniquely recover the curvature tensor from 
the function $(u,v)\mapsto R(u,v)(u+v)$, 
showing that the torsion of a connection,
and in case that is zero, the curvature of a connection, can be uniquely 
recovered from the knowledge of all the gaps.

If both the torsion tensor and the curvature tensor are identically
zero, then by Theorem \ref{frame flow bracket}, the 
horizontal distribution on $E$ defined by the connection 
is {\it involutive}, that is, closed under bracket. 
Hence by the Frobenius theorem, it follows that 
the connection is {\it integrable}, that is, 
$E$ locally admits horizontal sections, 
and therefore there are no further
obstructions for small sized geodesic quadrilaterals to close,
where the upper bound on $|s|$ depends on the starting triple $(P,u,v)$. 
However, if such an $(M,\nabla)$ is not geodesically complete, then 
there can still be a gap in a geodesic quadrilateral, {\it even if} the 
quadrilateral exists. For example, 
if we take $M$ to be the universal cover of $\R^2 - \{ (0,0)\}$,
with Riemannian metric pulled up from the Euclidean metric $dx^2+ dy^2$
below, then the pull-back of any square in $\R^2 - \{ (0,0)\}$
which has winding number $1$ around $(0,0)$ is an open geodesic 
quadrilateral in $M$. But if $\nabla$ has zero torsion, zero curvature and 
if moreover $M$ is geodesically complete, 
then there is not even a global obstruction, 
that is, the geodesic quadrilaterals will close for all 
$(P,u,v)$ and $s$. From this it can be seen that 
such a pair $(M,\nabla)$ will be isomorphic to a quotient of 
the affine space $\R^d$ with its natural constant connection.


\section{Spheres and hyperboloids}

In this section we directly verify the Theorem \ref{geodesic gap} for 
(i) the sphere $x^2+y^2+z^2 = r^2$ in the Euclidean space $\R^3$
with metric $dx^2+dy^2+dz^2$, 
which has a constant curvature $1/r^2$, and 
(ii) 
the hyperboloid $t^2-x^2-z^2 = r^2$, $t>0$ 
in the Minkowski space $\R^{1,2}$
with metric $-dt^2 + dx^2+ dy^2$,  
which has a constant curvature $-1/r^2$.

\stm {\bf Sphere in Euclidean $\R^3$.}
Let $\R^3$ be given the orientation in which the standard basis
$(e_1,e_2,e_3)$ is right-handed. 
We define $M \subset \R^3$ by the equation $x^2+y^2+z^2 =1$.
Let $M$ be given the orientation under which 
$(e_2,e_3)$ is a right-handed basis at $e_1\in M$.
The metric tensor on $M$ is induced from the Euclidean metric tensor
$dx^2+dy^2+dz^2$ on $\R^3$.
The data $(P,u,v)$ that consists of a point $P\in M$ and 
a right-handed orthonormal basis $(u,v)$ for $T_PM$ is the same 
as a right-handed orthonormal basis $(P,u,v)$ for $\R^3$, and 
so all such triples form the manifold $N$ which is 
a principal homogeneous space under the action of $SO(3)$.
Under this action, a group element $A\in SO(3)$ 
takes any triple  $(P,u,v)$ to the new triple $(AP,Au,Av)$.
We will identify $N$ with $SO(3)$ by representing each triple
$(P,u,v)$ by the matrix $X\in SO(3)$ whose columns are $P$, $u$ and $v$
respectively, so that in the parlance of the Introduction, the location
of $X$ is the first column of $X$. 

For $s\in \R$ and $X = (P,u,v)\in N$, recall that $\T_s(X)$
is the new triple $(P',u',v')$ as defined in the Introduction.
As the action of $SO(3)$ on $M$ preserves the metric, 
it preserves the geodesics and parallel transport of tangent vectors.
Hence $\T_s(X) = \T_s(XI) = X\T_s(I)$, which shows that the action 
of $\T_s$ on any $X$ is given by right multiplication by $A(s) = \T_s(I)$.

We now evaluate the matrix $A(s) = \T_s(I)$. The triple $(P_0,u,v)$ 
corresponding to $I$ has  
$P_0 = (1,0,0) = e_1,\, u = (0,1,0) = e_2,\, v = (0,0,1)=e_3$. 
The geodesics on $M$ are the great circles, and the distance along these 
is given by the angle subtended at the origin in $\R^3$.
Hence it is immediate that 
$$%
A(s) =
\left(\begin{array}{ccc}
\cos(s) & 0 & \sin(s)\\
\sin(s) & 0 & -\cos(s)\\
0          & 1 & 0
\end{array}
\right)
$$
Hence the vertices $P_i$ of the geodesic quadrilateral are the first columns 
of the powers $A(s)^i$, and the vertices $Q_i$ are
the first columns of the negative powers $A(s)^{-i}$.
This gives
$$P_2(s) = (\cos^2(s), \cos(s)\sin(s), \sin(s)),~~
Q_2(s) = (\cos^2(s), \sin(s), \cos(s)\sin(s)),$$
hence $Q_2(s) - P_2(s) = 
(0, \sin(s)(1-\cos(s)), - \sin(s)(1-\cos(s)))$. 
It follows that 
$$\lim_{s\to 0}\,{ Q_2(s) - P_2(s)\over s^3}
= (0, 1/2, -1/2) = {1\over 2}(u-v)$$
where $u,v\in T_PM$ are the vector $e_2,e_3$. This
proves the result for the gap between $P_2(s)$ and $Q_2(s)$. 

The first column of $A(s)^4$ is 
$$P_4 = \left[\begin{array}{c}
\cos^4(s) + 2\cos(s)\sin^2(s)\\
\cos^3(s) \sin(s) - \cos^2(s) \sin(s) + \sin^3(s) \\
\cos^2(s) \sin(s) - \cos(s) \sin(s)
\end{array}
\right]
$$
This implies that  
$$\lim_{s\to 0}\,{P_4(s) - P_0 \over s^3} = (0, 1/2, -1/2) = {1\over 2}(u-v).$$

Finally, to get the result for a sphere $x^2+y^2+z^2 = r^2$, which has 
constant Gaussian curvature $1/r^2$, we replace the variable $s$ by $s/r$
in $A(s)$, and multiply the first column of a power of $A(s/r)$ by $r$
to get the vertices $P_i,Q_i$. This gives
$$ \lim_{s\to 0}\,{P_4(s) - P_0 \over s^3} = 
\lim_{s\to 0}\,{ Q_2(s) - P_2(s)\over s^3}
= (0, 1/2r^2, -1/2r^2) = {\kappa\over 2}(u-v)$$

\stm {\bf Hyperboloid in Minkowskian $\R^{1,2}$.}
Let $\R^{1,2}$ be the Minkowski space with pseudo-Riemannian metric
$-dt^2 + dx^2 + dy^2$, and orientation chosen such that the basis
$e_t = (1,0,0), \,e_x = (0,1,0), \, e_y = (0,0,1)$ is right-handed.
Let $M\subset \R^{1,2}$ be defined 
by $t^2 - x^2 -y^2 = 1$ and $t>0$. The induced metric tensor 
$g$ on $M$ is positive definite, and has constant Gaussian curvature $-1$.
Let $M$ be given the orientation under which the basis $(e_x , e_y)$
of $T_{e_t}M$ is right-handed. 

Note that $M$ has a transitive free action of the proper orthochronous
Lorentz group $L^{\upa}_+$ which is the connected component of $SO(1,2)$.
Analogous to the previous example, all triples $(P,u,v)$, where $P\in M$
and $(u,v)\in T_PM$ is a right-handed orthonormal basis, form a
principal $L^{\upa}_+$-space. Again, such a triple $(P,u,v)$ can be identified
with the matrix $X$ in $L^{\upa}_+$ whose columns are $P,u,v$ respectively.
Thus, $N$ can be identified with $L^{\upa}_+$, and the action of $L^{\upa}_+$ 
becomes left multiplication. The action of $L^{\upa}_+$ preserves the metric, 
the geodesics, and parallel transport of vectors, as before.
So once again, we take $I$ as our base triple, and find that 
if we put $B(s) = \T_s(I)$, then for any triple $X$, we must have
$\T_s(X) = XB(s)$. 

As the geodesics and parallel transport have a simple description 
starting from the triple $I$, we can see that   
$$%
B(s) =
\left(\begin{array}{ccc}
\cosh(s)   & 0 & -\sinh(s)\\
\sinh(s)   & 0 & -\cosh(s)\\
0          & 1 & 0
\end{array}
\right)
$$
Hence starting from the triple $(P_0= e_t,u=e_x,v=e_y)$, 
the vertices $P_i$ of the geodesic quadrilateral are the first columns of
the powers $B(s)^i$, and the vertices $Q_i$ are
the first columns of the negative powers $B(s)^{-i}$.
This gives
$$P_2 = (\cosh^2 (s), \cosh (s) \sinh(s),  \sinh(s)),~ 
Q_2 = (\cosh^2 (s), \sinh(s),  \cosh (s) \sinh(s))$$
and
$$P_4 = \left[\begin{array}{c}
\cosh^4(s) - 2\cosh(s)\sinh^2(s)\\
\cosh^3(s)\sinh(s) - \cos^2(s) \sin(s) - \sin^3(s) \\
\cosh^2(s)\sinh(s) - \cosh(s)\sinh(s)
\end{array}
\right]
$$
Now a direct calculation gives the results
$$\lim_{s\to 0}\,{P_4(s) - P_0\over s^3} = 
\lim_{s\to 0}\,{ Q_2(s) - P_2(s)\over s^3} = -{1\over 2}(u-v).$$
A scaling by the factor $r$, where $B(s)$ is replaced by $B(s/r)$
and the points $P_i$ and $Q_i$ are replaced by $r$-times the 
first column of powers of $B(s/r)$
gives the result for the hyperboloid defined by
$t^2-x^2-y^2 = r^2$ in $\R^{1,2}$, which has curvature $-1/r^2$.

{\footnotesize 

\stm{\bf The case of an arbitrary surface metric $g$.} 
Let $P_0 \in M$, and let the Gaussian curvature of the given metric 
$g$ on $M$ take 
the value $\kappa_0 =\kappa(P_0)$ at $P_0$. 
There exists a coordinate neighbourhood $U$ of $P_0$ 
with local coordinates
$x,y$, known as Riemann normal coordinates centered at $P_0$,
such that $P_0 = (0,0)$, and the metric tensor takes the form
$$g = dx^2 + dy^2 -{\kappa_0\over 3}(y^2dx^2 - 2xydxdy +x^2dy^2) 
+ h_1(x,y)dx^2 + h_2(x,y)dxdy + h_3(x,y)dy^2$$
where each $h_i(x,y)$ is of order $\ge 3$ in $x$, $y$, that is,
$h_i(x,y) \in \mm^3 \subset C^{\infty}(U)$
where $\mm \subset C^{\infty}(U)$ is the maximal ideal generated by $x,y$,
(which is the set of all $C^{\infty}$-functions on $U$
that vanish at $P_0$). 
Hence, up to terms of $2$nd order (means modulo $\mm^3$),
the pair $(U,g)$ is the same as the pair $(U,g')$ where $g'$ is 
a Riemannian metric on $U$ with constant Gaussian curvature $\kappa_0$, as
in Riemann normal coordinates, both $g$ and $g'$ have the common form
$$g = dx^2 + dy^2 -{\kappa_0\over 3}(y^2dx^2 - 2xydxdy +x^2dx^2) 
\mbox{ mod }\mm^3 Sym^2(TM)$$
One can heuristically argue from this that the 
points $P_i(s)$, $Q_i(s)$ for the two metrics $g$ and $g'$ are congruent 
modulo $s^4$, and consequently the quadrilateral
gaps $P_4(s) - P_0$ or $Q_2(s) - P_2(s)$ for the two metrics 
are congruent modulo $s^4$. Hence the validity of the gap formula
$P_4(s) - P_0 = Q_2(s) - P_2(s) = s^3(\kappa_0/2)(u-v) \mbox{ mod }s^4$
for constant curvature metrics, which we verified above, 
implies its validity for all Riemannian surfaces.
We do not make this argument rigorous here, instead, we deduce the 
validity of the gap formula for arbitrary Riemannian surfaces from
the Theorem \ref{geodesic gap}, as explained in the Introduction.
 }

\rm

\break

\section{Brackets, flows and gaps}

\centerline{\it Taylor expansion and equivalence modulo $s^n$}

We recall some elementary facts about rings of real valued smooth 
functions. All smooth functions on $(-a,a)$ that vanish 
at $s =0$ form the principal ideal $(s) = sC^{\infty}(-a,a)$ in the ring
$C^{\infty}(-a,a)$, and all smooth functions $f(s,t)$ 
on $D = (-a,a)\times (-a,a)$
that vanish at $(0,0)$ form the ideal $(s,t)\subset C^{\infty}(D)$.
If $W$ is any open subset of some $\R^n$ 
(or more generally, a smooth manifold), then all smooth 
$f: (-a,a)\times W \to \R$ that vanish on $\{ 0 \} \times W$ 
form the principal ideal $(s) = sC^{\infty}((-a,a)\times W)$ in the ring
$C^{\infty}((-a,a)\times W)$, and all smooth functions $f$ 
on $D\times W$
that vanish at $\{ (0,0)\} \times W$ form the ideal 
$(s,t)\subset C^{\infty}(D\times W)$.
The following lemma lists some frequently used elementary 
facts.


\lemma\label{Taylor}{\it
(1) Taylor expansion: For any 
smooth $f$ on $(-a,a)\times W$ and any $n\ge 0$, there exist unique functions
$g_0,\ldots,g_n\in C^{\infty}(W)$ such that 
$f = g_0+ g_1s+ \ldots + g_ns^n \mbox{ mod }(s^{n+1})$ in the ring 
$C^{\infty}((-a,a)\times W)$, where
$g_0(w) = f(0,w)$ and $g_i(w) =  {1\over i!}{\pa^i f\over \pa 
s^i}(0,w)$ for $i = 1,\ldots, n$. 

(2) In particular, if for each $w_0\in W$ we have $f(s,w_0)\in 
s^nC^{\infty}(-a,a)$, then $f \in s^nC^{\infty}((-a,a)\times W)$.

(3) Both the above statements hold with $(-a,a)$ replaced by 
$D$, where we have the Taylor expansion
$$f = \sum_{0\le j+k\le n} h_{jk}\, s^jt^k \mbox{ mod }(s,t)^{n+1}$$  
in the ring $C^{\infty}(D\times W)$, where 
$$h_{jk}(w) = {1\over j!k!}{\pa^{j+k} f \over \pa s^j\pa t^k}(0,0,w)
\in C^{\infty}(W).$$

(4) For $a,b >0$, all functions on $U = (-a,a)\times (-b,b)$ that vanish on the 
diagonal $\Delta \subset U$ (which is defined by $s=t$) form the principal 
ideal $(s-t)$ in $C^{\infty}(U)$. A similar statement holds for 
$\Delta \times W \subset U \times W$.
}

\rm

\bigskip

\centerline{{\it Brackets and gaps in flows}}

For vector fields $\xi$ and $\eta$ on a differential manifold $M$, 
the bracket $[\xi,\eta] = \xi\circ \eta - \eta\circ \xi$
can be seen in terms of two different kinds of quadrilateral gaps
that are related to their flows. Let $\varphi_s$ and $\psi_s$ be
the respective flows, defined locally for small values of $s$.
In particular, if $U\subset M$ is an open subset whose closure is compact,
then there exists some $a>0$ such that both $\varphi_s$ and $\psi_s$ 
are defined as smooth maps $(-a,a)\times U \to M$.

If $P_0\in M$ and if 
$s$ is small enough, the following construction
makes sense. Let $P_1(s) = 
\varphi_s(P_0)$, $P_2(s) = \psi_s(P_1(s))$, 
$P_3(s) = \varphi_{-s}(P_2(s))$ and 
$P_4(s)= \psi_{-s}(P_3(s))$.
Let $Q_1(s) = \psi_s(P_0)$ and $Q_2(s) = \varphi_s(Q_1(s))$. 
These points form two open quadrilaterals.

The first open quadrilateral 
has the five vertices $P_0,P_1(s),\ldots, P_4(s)$, with
successive vertices joined by integral curves of $\xi$,
$\eta$, $-\xi$ and $-\eta$, respectively.
The second open quadrilateral has the five successive   
vertices $Q_2(s), Q_1(s), P_0, P_1(s), P_2(s)$, with 
successive vertices joined by integral curves of $-\xi$, $-\eta$,
$\xi$ and $\eta$, respectively.
As in general the two flows may not commute, 
$P_4(s) \ne P_0$ and $P_2(s) \ne Q_2(s)$ in general. 

{
\begin{center}
\includegraphics[scale=.15]{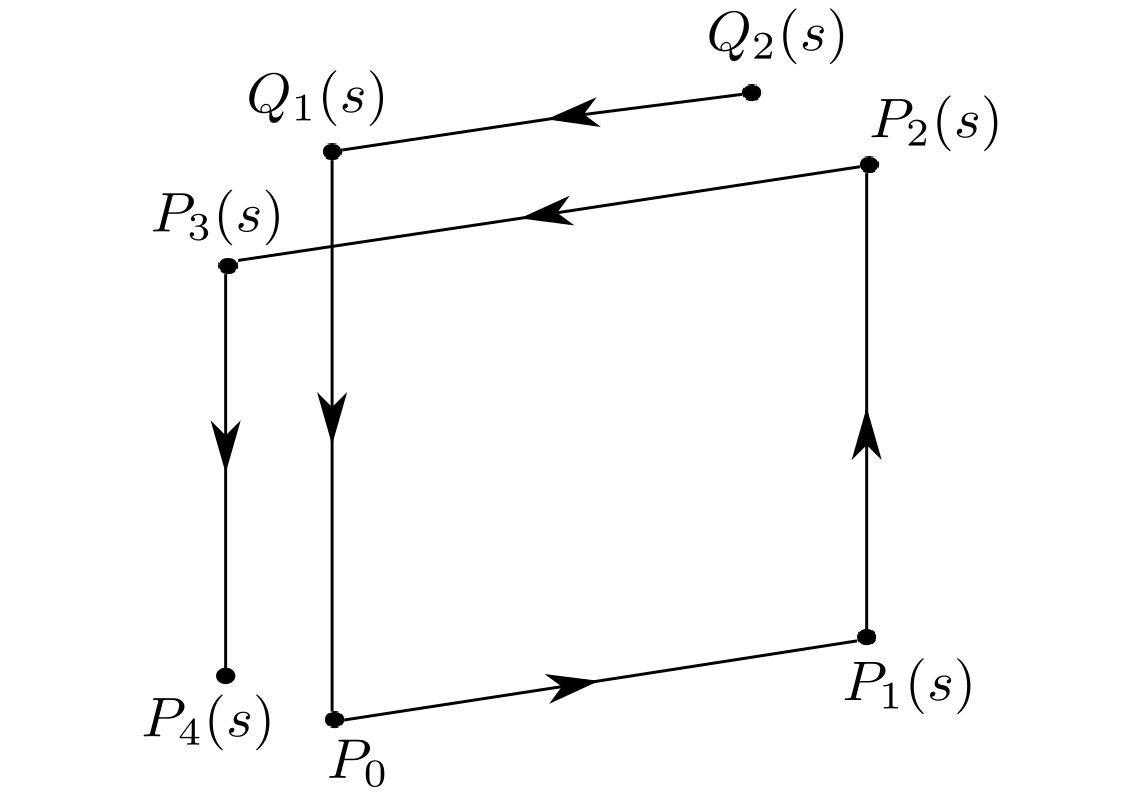}
\end{center}
}

Given any point $x\in M$, there exists a real number $a>0$ and 
a neighbourhood $W$ of $x$ in $M$, such that for any $s\in (-a,a)$ and 
$P_0\in W$, the corresponding 
points $P_1(s), P_2(s), P_3(s),P_4(s),Q_1(s),Q_2(s)$
are defined and are smooth functions from $(-a,a)\times W$ to $M$.
This is a consequence of the basic existence and uniqueness theorem for 
ordinary differential equations.

Given any smooth function $f$ in a neighbourhood of $x\in M$, we now
define the following two functions on $(-a,a)\times W$
where $a>0$ and $W$ is a neighbourhood of $x$ in $M$, where
both $a$ and $W$ are chosen to
be sufficiently small. 
The {\bf flow gap functions} 
$F_I(P_0,\xi,\eta,f,s)$ and $F_{II}(P_0,\xi,\eta,f,s)$ are the
functions from $(-a,a)\times W$ to $M$ defined by 
$$F_I(P_0,\xi,\eta,f,s) = f(P_4(s)) - f(P_0),\mbox{ and }
  F_{II}(P_0, \xi,\eta,f, s) = f(P_2(s)) - f(Q_2(s)).$$

We now rephrase the Lemma \ref{meaning of bracket} 
in terms that will be more useful to us.

\lemma\label{meaning of bracket again}
{\it
With notation as above, we have the following relations 
in the ring $C^{\infty}((-a,a)\times W)$ of all smooth real functions 
on $(-a,a)\times W$, where $(s^3)$ is the ideal consisting of all
multiples of $s^3$.\\
(i) $F_I(P_0,\xi, \eta,f,s) = s^2[\xi,\eta]_{P_0}(f) \mbox{ mod }(s^3)$, and \\
(ii) $F_{II}(P_0,\xi, \eta,f,s) = s^2[\xi,\eta]_{P_0}(f) \mbox{ mod }(s^3)$.
}
\rm

{\bf Strategy of the proof of the Lemma \ref{meaning of bracket again}.} 
To begin with, note that as a consequence of the Lemma \ref{Taylor}.(2),
it is enough to prove the above lemma for a fixed $P_0$, instead of 
allowing $P_0$ to vary over an open set $W$. 
Though the limits corresponding to them
are equal, each of the two gap functions
$F_I(P_0,\xi,\eta,f,s)$ and $F_{II}(P_0,\xi,\eta,f,s)$ has a different advantage.
Using $F_I$, we see that as $s$ varies the point $P_4(s)$  
moves on a smooth curve through $P_0$, whose tangent at $P_0$, when 
calculated w.r.t. the new parameter $s^2$, is $[\xi,\eta]_{P_0}$.
On the other hand, the definition of $F_{II}$  has more symmetry
and has fewer compositions of flows, so $F_{II}$ is much easier to
calculate. Given this, we will first prove the statement about $F_{II}$
in Lemma \ref{meaning of bracket again}. Though a direct proof of the statement 
about $F_I$ is possible along the same lines by working a bit harder, we
will instead use a different trick. This involves deducing
the statement (i) about $F_I$ for a fixed $P_0$ 
from the full form of the statement (ii) about $F_{II}$
in which $P_0$ is not fixed but varies over an open set $W\subset M$. 
A similar trick will be deployed later in the proof of the 
Theorem \ref{geodesic gap}, where we first prove the statement about the 
gap function $G_{II}$ related to the geodesic quadrilateral
$Q_2,Q_1,P_0,P_1,P_2$, and then use this trick to deduce the statement
about the gap function $G_I$ related to the geodesic quadrilateral
$P_0,P_1,P_2,P_3,P_4$, 
which is considerably more difficult (being 
much longer and computationally messy) 
to approach directly.

{\bf Proof of Lemma \ref{meaning of bracket again}.(ii):} 
As explained above, it is enough to prove the statement for 
a fixed $P_0$, instead of allowing it to vary over $W$.
For small enough $s,t$ let 
$$\P(s,t) = \psi_t\circ\varphi_s (P_0) \in M,$$
in particular, $P_0 = \P(0,0)$, $P_1(s) = \P(s,0)$ and $P_2(s) = \P(s,s)$. 
Let $f$ be a smooth function on a neighbourhood of $P_0$. 
Then regarded as a function on $(-a,a)\times (-a,a)$, we have
$$
{\pa\over \pa s}f(\P(s,0)) = \xi_{\P(s,0)}(f), \mbox{ and }
{\pa\over \pa t}f(\P(s,t)) = \eta_{\P(s,t)}(f)$$
By expansions in $s$ and $t$, we have the following equalities. 
\begin{eqnarray*}
f(\P(s,0)) & = &   f(P_0) + s \xi_{P_0}(f)  + {s^2\over 2}\xi_{P_0}\xi(f) 
+ s^3h_1(s),\\
f(\P(s,t)) & = &  f(\P(s,0)) + t \eta_{\P(s,0)}(f) + 
                {t^2\over 2}\eta_{\P(s,0)}(\eta(f)) + t^3 h_2(s,t),\\
\eta_{\P(s,0)}(f) & = & \eta_{P_0}(f) + s \xi_{P_0}(\eta(f)) + s^2h_3(s),\\
\eta_{\P(s,0)}(\eta(f)) & = & \eta_{P_0}(\eta(f)) + s h_4(s)
\end{eqnarray*}
where $h_1(s)$, $h_3(s)$ and $h_4(s)$ are smooth real functions
on $(-a,a)$, while $h_2(s,t)$ is a  smooth real function
on $(-a,a)\times (-a,a)$.
Hence we get 
$$f(\P(s,t)) - f(P_0) = s \xi_{P_0}(f) + t \eta_{P_0}(f) +  
st (\xi_{P_0}\eta(f)) 
+ {s^2\over 2}\xi_{P_0}\xi(f) +{t^2\over 2} \eta_{P_0}\eta(f) 
+ H(s,t)$$
where the error term $H(s,t)$ is given by
$$H(s,t) = s^3h_1(s) + t^3 h_2(s,t) + s^2t h_3(s) + {st^2\over 2} h_4(s).$$
Now put $t=s$ in the above, that is, restrict the above function 
to the diagonal $\Delta  \subset (-a,a)\times (-a,a)$.
As 
$P_1(s)= \P(s,0)$ and $P_2(s)= \P(s,s)$, we get 
$$f(P_2(s)) - f(P_0) = s(\xi_{P_0}(f) + \eta_{P_0}(f)) +  
s^2 (\xi_{P_0}\eta(f)) 
+ {s^2\over 2}(\xi_{P_0}\xi(f) + \eta_{P_0}\eta(f))
+ s^3h(s)$$
where $h(s)$ is a smooth function on $(-a,a)$.
Similarly, by first travelling along the flow of $\eta$ to 
reach $Q_1$ and then along the flow of $\xi$ to reach $Q_2$ 
(that is, by just interchanging $\xi$ and $\eta$ in the above expression), 
we get
$$f(Q_2(s)) - f(P_0) = s (\eta_{P_0}(f) + \xi_{P_0}(f)) +  
s^2(\eta_{P_0}\xi(f)) 
+ {s^2\over 2}(\eta_{P_0}\eta(f)+ \xi_{P_0}\xi(f))  
 + s^3k(s) $$
for some smooth function $k(s)$ on $(-a,a)$.
Taking the difference, 
$$f(P_2(s)) - f(Q_2(s)) 
= s^2(\xi_{P_0}\eta(f) - \eta_{P_0}\xi(f)) \mbox{ mod }(s^3)
= 
s^2[\xi, \eta]_{P_0}(f) \mbox{ mod }(s^3),$$
which proves the statement of Lemma \ref{meaning of bracket again}.(ii).

{\bf The trick which gets $F_I$ from $F_{II}$.}
We have proved that for any $x\in M$, there is a neighbourhood
$x\in W\subset M$ and an $a>0$ such that both the quadrilaterals are 
defined for any $P_0\in W$ and $|s|<a$, and we have 
$$F_{II}(P_0,\xi, \eta,f,s) = s^2[\xi,\eta]_{P_0}(f) \mbox{ mod }
s^3C^{\infty}((-a,a)\times W),$$ 
which is the statement (ii) of the Lemma. 
To deduce from this the statement (i) of the Lemma, we begin by noting 
the identity
$$F_I(P_0,\xi,\eta,f,s) = F_{II}(P_2(s),-\xi,-\eta,f,s)$$
which is immediate from the definitions. 
Now for a fixed $P_0$, consider the map $(-a,a)\to M$
that sends 
$\tau \mapsto 
P_2(\tau) = \psi_{\tau}\circ\varphi_{\tau}(P_0)$. 
As $P_2(0) = P_0$, by continuity there exists $0< b < a$
such that $P_2(\tau) \in W$ whenever $\tau \in (-b,b)$.  
The smooth map 
$$(-b,b)  \times (-b,b) \to (-a,a)\times W : 
(s, \tau)\mapsto (s, P_2(\tau))$$
induces a ring homomorphism 
$C^{\infty}((-a,a)\times W) \to C^{\infty}((-b,b) \times (b,b))$
under which the 
relationship in the ring $C^{\infty}((-a,a)\times W)$ in the statement of
the Lemma \ref{meaning of bracket again}.(ii) 
gives the following relationship in the ring $C^{\infty}((-b,b) \times (-b,b))$:
$$F_{II}(P_2(\tau),\xi,\eta,f,s) = s^2[\xi,\eta]_{P_2(\tau)}(f) \mbox{ mod }
s^3C^{\infty}((-b,b) \times (-b,b)).$$
It should be noted that the above step used the full form of 
Lemma \ref{meaning of bracket again}.(ii), in which $P_0$ varies
over a neighbourhood $W$ of $x\in M$, as we needed to apply it to the 
variable base point $P_2(\tau)\in W$.

Under the ring homomorphism $C^{\infty}((-b,b) \times (-b,b))\to C^{\infty}(-b,b)$ 
induced by the inclusion of the diagonal, under which  
$s\mapsto s$ and $\tau \mapsto s$, the extension of the ideal 
$s^3 C^{\infty}(-b,b) \times (-b,b))$ is the ideal 
$s^3C^{\infty}(-b,b)$. 
Applying this homomorphism
to the above equality (means, putting $s=\tau$), we get 
$$F_{II}(P_2(s),\xi,\eta,f,s) = s^2[\xi,\eta]_{P_2(s)}(f) \mbox{ mod }
s^3C^{\infty}(-b,b).$$
Now, $P_2(s)$ is a smooth function of $s$, and for $s=0$ it takes 
the value $P_0$. Hence by Lemma \ref{Taylor}.(1) we have 
$$[\xi,\eta]_{P_2(s)}(f) =  [\xi,\eta]_{P_0}(f) \mbox{ mod }
sC^{\infty}(-b,b).$$
Hence substitution gives us 
$$F_{II}(P_2(s),\xi,\eta,f,s) = s^2[\xi,\eta]_{P_0}(f) \mbox{ mod }
s^3C^{\infty}(-b,b).$$

Hence we finally have
\begin{eqnarray*}
F_I(P_0,\xi,\eta,f,s)  
& =& F_{II}(P_2(s),-\xi,-\eta,f,s)\\
& = & s^2[-\xi,-\eta,]_{P_0}(f) \mbox{ mod }s^3C^{\infty}(-b,b)\\
& = & s^2[\xi,\eta,]_{P_0}(f) \mbox{ mod }s^3C^{\infty}(-b,b) 
\mbox{ as desired.} 
\end{eqnarray*}


\section{Connections: basic notions}
In this section we first recall 
some well known basic notions related to affine connections (see 
e.g. [1] for a more comprehensive introduction), and then formulate 
the Lemma \ref{curvature identities}.

An {\bf affine connection} on a smooth manifold $M$ associates to
any tangent vector fields $\xi$ and $\eta$ on an open subset $U\subset M$
a new tangent vector field $\nabla_{\xi}\eta$ on $U$. 
This is $C^{\infty}(U)$-linear in the variable $\xi$, while it is $\R$-linear
in $\eta$ and satisfies the Leibniz rule: $\nabla_{\xi}(f\eta)
= \xi(f)\eta + f\nabla_{\xi}(\eta)$ for all $f\in C^{\infty}(U)$. In particular,
if $P\in U$ and $u\in T_PM$, then the vector $\nabla_u\eta \in T_PM$, 
which is defined to be $(\nabla_{\xi}\eta)_P$ where $\xi$ is any tangent vector 
field in an open neighbourhood of $P$ such that $\xi_P = u$, is well defined.

Let $M$ be a manifold with an affine connection $\nabla$. 
Let $\gamma(s)$ be a parameterized curve in $M$ defined on an interval
$s\in (-a,a)\subset \R$ (this means $\gamma: (-a,a)\to M$ is smooth). 
A smooth vector field $v(s)$ along $\gamma$ consists of giving
a vector $v(s) \in T_{\gamma(s)}M$ for all $s\in (-a,a)$,
which is smooth as a function of $s$. For example, the {\bf tangent vector
field of $\gamma$}, denoted by $\dot{\gamma}$ or $d\gamma/ds$, 
is a smooth vector field 
along $\gamma$. If $\nabla$ is an affine connection on $M$, 
and $\gamma$ and $v$ are as above, then we can define 
the {\bf covariant derivative} of $v$
along $\gamma$ to be a certain smooth tangent vector field along $\gamma$,
denoted by $\nabla_{\dot{\gamma}}v$ or $\nabla_{d\gamma/ds}v$. 
While this can be naturally defined in a coordinate-free manner by
working over the graph of $\gamma$ in $(-a,a)\times M$, what we will 
find more useful is its local coordinate expression, which we next
describe.

For this, let $U$ be a coordinate chart around $P$ with coordinates
$(x^1,\ldots, x^d)$. 
The {\bf connection coefficients} of $\nabla$ are the real-valued smooth 
functions $\Gamma^i_{jk}$ on $U$ defined by
$$\nabla_{\pa\over \pa x^k} (
{\textstyle{\pa\over \pa x^j}}) = \Gamma^i_{jk}{\textstyle{\pa\over \pa x^i}}.$$
where we have used the {\bf summation convention} under which 
there is understood to be a summation over an index which is repeated with
one occurance as a subscript and another as a superscript. The superscript
$i$ on $x^i$ is to be regarded as a subscript in the expression $\pa/\pa x^i$. 
The summation signs are suppressed.

Let the point $\gamma(s)$ have coordinates $(x^1(s),\ldots, x^d(s))$, so that
the tangent vector is $\dot{\gamma} = {dx^i\over ds}{\pa\over \pa x^i}$. 
Let $v(s) = v^i(s){\pa\over \pa x^i}$. 
Let $\Gamma^i_{jk}(s) = \Gamma^i_{jk}(\gamma(s))$. 
Then 
$$(\nabla_{\dot{\gamma}}v)(s)  = \left(
{d v^i\over ds}(s) + \Gamma^i_{jk}(s)v^j(s)
{dx^k\over ds}(s) \right)
{\pa\over \pa x^i}|_{\gamma(s)},$$
which we write more briefly as
$(\nabla_{\dot{\gamma}}v)^i = {d v^i\over ds} + \Gamma^i_{jk}v^j\dot{\gamma}^k$.
We say that the vector field $v$ along $\gamma$ is {\bf parallel transported} 
along $\gamma$ if $\nabla_{\dot{\gamma}}v=0$, that is,
${d v^i\over ds}  = - \Gamma^i_{jk}v^j\dot{\gamma}^k$.
The parameterized curve $\gamma(s)$ is called a {\bf geodesic}
if its tangent vector field $\dot{\gamma}$ is   parallel transported
along $\gamma$, that is, 
$\nabla_{\dot{\gamma}}\dot{\gamma} = 0$. In coordinate terms, 
this is the equation
$${d^2x^i\over ds^2} = - \Gamma^i_{jk}{dx^j\over ds}{dx^k\over ds}.$$
Given any $P\in M$ and $u\in T_PM$, for a small enough $a>0$ there is a 
unique geodesic $\gamma: (-a,a) \to M$ such that 
$\gamma(0) = P$ and ${d\gamma/ds}(0) = u$.

\bigskip

\centerline{\it Coordinates and vector fields on the frame bundle}

We will denote by $\pi: E\to M$ the frame bundle associated to $TM$, 
which is a principal $GL(d)$-bundle where $d = \dim M$,
whose fiber over $x\in M$ is the set of all linear bases for $T_xM$. 
If $(x^1, \ldots,x^d)$ are local coordinates
on $M$ defined on an open subset $U\subset M$, then on the open
subset $\pi^{-1}(U)\subset E$, we get coordinates
$(x^1, \ldots,x^d, \x^1_1, \ldots,\x^i_j,\ldots, \x^d_d)$. 
These are defined as follows. If $x\in U$ has coordinates
$(a^1, \ldots,a^d)$, and 
if $y=(u_i,\ldots,u_d)$ is a basis
for $T_xM$, with $u_j = b^i_j\pa_i$ where $\pa_i = \pa/\pa x^i$, 
then $\pi(y) = x$, and the coordinates of $y$ are 
$x^i = a^i$ and $\x^i_j = b^i_j$.

The group $GL_d(\R)$ acts on the right on $E$, under which 
a matrix $A = (A^p_q) \in GL_d(\R)$ moves the point $y = (u_1,\ldots,u_d)\in E$
to the point $(v_1,\ldots,v_d) = (u_1,\ldots,u_d)A \in E$ where  
$v_i = u_jA^j_i$. A connection $\nabla$ on $TM$ is the same as a distribution
(vector subbundle) $D\subset TE$ which is (i) preserved by the action of 
$GL_d(\R)$, and (ii) is `horizontal', that is, supplementary to 
the kernel of $\pi_* : TE \to \pi^*TM$. As $\pi$ is a submersion, 
the kernel of $\pi_*$ is a vector subbundle of rank $d^2$. This is called
the `vertical subbundle', and it is naturally ismorphic to 
$\pi^*\un{End}(TM)$. 

We now define $d$ global vector fields $\xi_1, \ldots, \xi_d$ on 
$E$, which are everywhere linearly independent and span $D$. 
At any point $y = (u_1, \ldots, u_d) \in E$ over $x\in M$, the linear map 
$\pi_*(y): D_y \to T_xM$ is an isomorphism, so there is a 
unique $\xi_i \in D_y$ such that $$\pi_*(y) \xi_i = u_i$$
for all $i = 1, \ldots, d$.
In coordinate terms, if $\Gamma^i_{jk}$ are the connection coefficients for
the connection $\nabla$ on $TM$ in the local coordinates $(x^i)$, then 
$D_y$ has the basis $\xi_1, \ldots, \xi_d$ where
$$\xi_m = 
\x^i_m {\pa\over \pa x^i} 
- \Gamma^i_{jk}\x^j_{\ell}\x^k_m {\pa \over \pa \x^i_{\ell}}$$
In the above, there is a summation over the
indices $i,j,k,\ell$, which are repeated with 
one occurance as a subscript and another as a superscript. For this
purpose, the index $\ell$ in ${\pa \over \pa \x^i_{\ell}}$ is regarded
as a superscript. The summation sign is suppressed. 

\bigskip

\centerline{{\it How to recover the Riemann curvature tensor
from the data $R(u,v)(u+v)$}}

The {\bf torsion tensor} of an affine connection $\nabla$ associates
to any pair $(\xi,\eta)$ of tangent vector fields on an open subset 
$U\subset M$ the tangent vector field $T(\xi,\eta) = 
\nabla_{\xi}\eta - \nabla_{\eta}\xi - [\xi,\eta]$ on $U$. This turns
out to be $C^{\infty}(U)$-linear in both $\xi$ and $\eta$, justifying the name
`tensor'. In particular, if $P\in M$ and $u,v\in T_PM$, we get
a well-defined vector $T(u,v)\in T_PM$, which is $\R$-linear in each 
of $u,v$.
An affine connection is said to be {\bf symmetric} if the tensor $T$ 
is identically zero on $M$.

The {\bf Riemann curvature tensor} of $\nabla$ is the operator
that associates
to any triple $(\xi,\eta,\zeta)$ of tangent vector fields on an open subset 
$U\subset M$  the tangent vector field 
$R(\xi, \eta)\zeta = \nabla_{\xi}\nabla_{\eta}\zeta- \nabla_{\eta}\nabla_{\xi} \zeta
- \nabla_{[\xi,\eta]}\zeta$ on $U$. This turns
out to be $C^{\infty}(U)$-linear in each of the three variables, 
justifying the name `tensor'. In particular, if $P\in M$ and 
$u,v,w\in T_PM$, then we get a well-defined vector
$R(u,v)w \in T_PM$, which is $\R$-linear in each of $u,v,w$.

The curvature tensor satisfies 
$R(u,v)w = - R(v,u)w$ (skew-symmetry) and moreover 
if $\nabla$ is symmetric, then 
$R(u,v)w + R(v,w)u + R(w,u)v =  0$ (the algebraic Bianchi identity).
From these, the following lemma is immediate.

\lemma\label{curvature identities} 
{\it If $\nabla$ is symmetric, then 
for any three vectors $u,v,w \in T_P M$,
the following identities hold.
\begin{eqnarray*}
R(u,v)u & = & {\textstyle {1\over 2}}R(2u,v)(2u+v) - R(u,v)(u+v), 
                \mbox{ and}\\
R(u,v)w & = & {\textstyle {1\over 3}}(R(u, v+w)(v+w) - R(u,v)v - R(u,w)w \\
        &   &  + R(u+w,v)(u+w) - R(u,v)u - R(w,v)w).
\end{eqnarray*}
Hence  
the curvature operator $TM\times_M TM\times_M TM \to TM : (u, v, w) 
\mapsto R(u,v)w$ can be recovered uniquely from the 
map $TM\times_M TM \to TM : (u, v) \mapsto R(u,v)(u+v)$.
}

\rm

\section{Frame flow, torsion and curvature}

The vector fields $\xi_1,\ldots, \xi_d$ are so defined that the integral curve 
of $\xi_m$ on $E$ through 
a point $y = (u_1, \ldots, u_d) \in E$
projects under $\pi : E \to M$ to the geodesic
through $x = \pi(y)$ with tangent $u_m\in T_xM$, 
and the integral curve in $E$ 
parallel translates each of $u_1, \ldots, u_d$
along this geodesic through $x$. Hence we call $\xi_m$ as the
{\bf $m$th frame flow field}, and its flow on $E$ as the 
{\bf $m$th frame flow}. This is a natural generalization of the notion of
the geodesic flow on $TM$, under which only the tangent vector is 
parallel translated, instead of a full frame.

{}


\medskip

{\bf Proof of Theorem \ref{frame flow bracket}.} 
We will denote $\pa f/\pa x^i$ by $f_{,i}$.
We will use the notation
$$\pa_i = {\pa\over \pa x^i} \mbox{ and }D^p_q = {\pa\over \pa \x^q_p}$$
With this, we have 
$$\xi_m = \x^i_m\pa_i - \Gamma^i_{jk}\x^j_a\x^k_mD^a_i 
\mbox{ and } 
\xi_n = \x^p_n\pa_p - \Gamma^p_{qr}\x^q_s\x^r_n D^s_p$$
Note that $\pa_iD^p_q = D^p_q\pa_i$, $\pa_i(\x^p_q)=0$  
and $D^i_j(\x^p_q) = \delta^i_q\delta^p_j$
in terms of Kronecker symbols.
Hence we have 
\begin{align*}
\xi_m\circ \xi_n
=\ & (\x^i_m \pa_i - \Gamma^i_{jk}\x^j_a\x^k_m D^a_i)
(\x^p_n \pa_p - \Gamma^p_{qr}\x^q_s\x^r_n D^s_p)  \\
= \ & - \Gamma^i_{cb}\x^b_m \x^c_n \pa_i 
+ (- \Gamma^p_{q c,b}  + \Gamma^i_{qb}\Gamma^p_{ic}
+ \Gamma^i_{cb}\Gamma^p_{q i})\x^q_s\x^b_m\x^c_nD^s_p \\
& + \mbox{ second order partial derivative terms}. 
\end{align*}
Similarly, by exchanging $m$ and $n$ and renaming the dummy variables,
\begin{align*}
\xi_n\circ \xi_m
= \ & - \Gamma^i_{bc}\x
^b_m \x^c_n \pa_i +
(-\Gamma^p_{q b,c} + \Gamma^i_{qc}\Gamma^p_{ib}+ \Gamma^i_{bc}\Gamma^p_{q i})
   \x^q_s \x^b_m \x^c_nD^s_p \\
& + \mbox{ second order partial derivative terms}. 
\end{align*}
As the bracket is a differential operator of first order 
(which is because the second order terms
exactly cancel as partial derivatives commute), we get
\begin{align*}
\left[\xi_m, \xi_n \right] 
= \ & \xi_m\circ \xi_n - \xi_n\circ \xi_m \\
= \ &(- \Gamma^i_{cb}+ \Gamma^i_{bc})\x^b_m \x^c_n\pa_i  + 
\mbox{ linear combinations of the } D^s_p.
\end{align*}
Recall that at any point $y = (u_1,\ldots,u_n)\in E$ over $x\in M$,
we have $\pi_*(y)(\xi_r(y)) = u_r$, where $\pi_*(y) : T_yE \to T_xM$ 
is the derivative of $\pi: E \to M$. 
$T(\pa_m, \pa_n) = (\Gamma^i_{nm} - \Gamma^i_{mn})\pa_i$, that is, 
$T^i_{mn} = \Gamma^i_{nm} - \Gamma^i_{mn}$.
So taking $u_r = \pa_r$, which corresponds to $\x^s_r = \delta^s_r$,
the above gives
$$\pi_*[\xi_m,\xi_n] = - T^i_{mn}\pa_i = - T(\pa_m,\pa_n)$$ 
As the $\pa_r$ form a basis for $T_xM$, this shows that 
the torsion of $\nabla$ is given by 
$T(u_i,u_j) = - \pi_*(y) [\xi_i,\xi_j]$ as claimed. 
This proves the statement (1) in the Theorem.

Now assuming $\nabla$ is symmetric, that is, $\Gamma^i_{jk} = \Gamma^i_{kj}$, 
we get 
\begin{align*}
\left[\xi_m, \xi_n \right] 
= \ & \xi_m\circ \xi_n - \xi_n\circ \xi_m \\
= \ &  ((- \Gamma^p_{q c,b}  + \Gamma^i_{qb}\Gamma^p_{ic}
+ \Gamma^i_{cb}\Gamma^p_{q i}) -  
(-\Gamma^p_{q b,c} + \Gamma^i_{qc}\Gamma^p_{ib}+ \Gamma^i_{bc}\Gamma^p_{q i}))
\x^q_s \x^b_m \x^c_nD^s_p \\
= \ &  
(\Gamma^p_{q b,c}- \Gamma^p_{q c,b} + \Gamma^i_{qb}\Gamma^p_{ic}
- \Gamma^i_{qc}\Gamma^p_{ib}) \x^q_s \x^b_m \x^c_nD^s_p\\
=  \ & - R^p_{qbc}\x^q_s \x^b_m \x^c_nD^s_p 
\end{align*}
where 
$$R^p_{qbc} = \Gamma^p_{q c,b} - \Gamma^p_{q b,c} 
+ \Gamma^i_{qc}\Gamma^p_{ib} - \Gamma^i_{qb}\Gamma^p_{ic}
$$ 
are the coordinates of the Riemann curvature tensor when $\nabla$ is symmetric.
Taking $\x^j_i = u_i^j$, it follows that 
$$R(u_i,u_j)_x = - [\xi_i,\xi_j]_y$$
as elements of $End(T_xM)$. 
This completes the proof of (3), and hence that of the 
Theorem \ref{frame flow bracket}.

\section{Quadrilateral gap for an affine connection}

Let there be given 
tangent vectors $u_0,v_0\in T_{P_0}M$ at some point $P_0\in M$. 
For small values of $s$, recall from the Introduction that by repeated
application of the operator $\T_s$ or its inverse to the triple
$(P_0,u_0,v_0)$ we defined two open geodesic
quadrilaterals $P_0,P_1(s),P_2(s),P_3(s),P_4(s)$ 
and $Q_2(s), Q_1(s), P_0, P_1(s), P_2(s)$. 
The successive vertices in each quadrilateral 
are joined by geodesic segments as described earlier.

{
\begin{center}
\includegraphics[scale=.25]{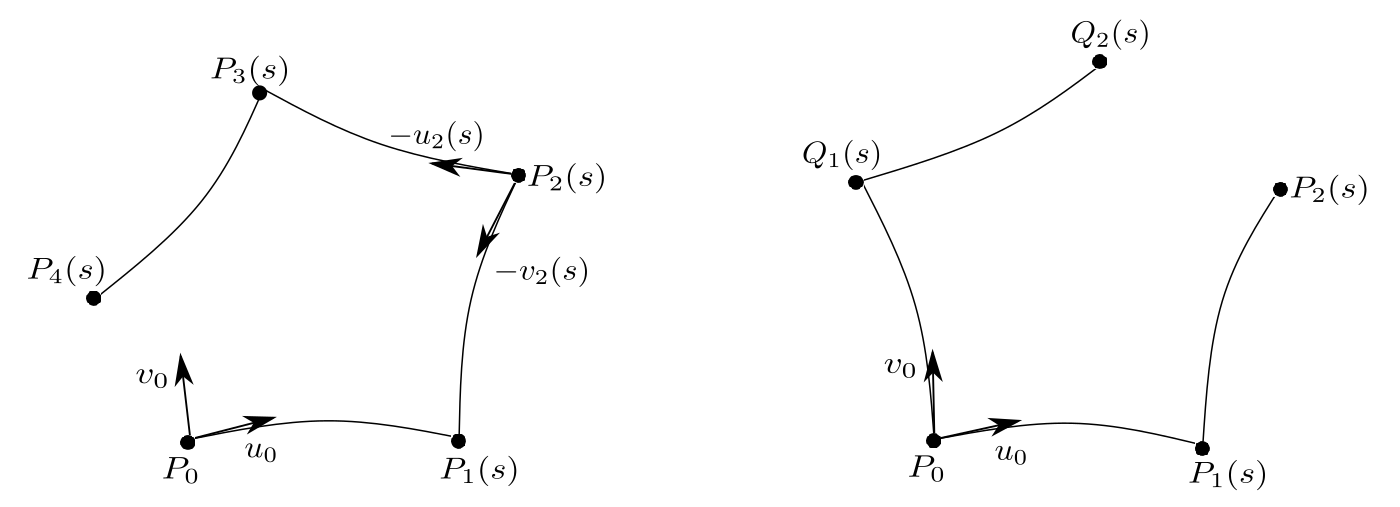}
\end{center}
}

By the fundamental existence and uniqueness 
theorem for ODEs with smooth parameters, applied
to the geodesic equation and the parallel transport equation, 
for $(x,u_x,v_x) \in TM\times_MTM$ there exists 
an open neighbourhood $W\subset TM\times_MTM$ and a real number $a>0$ 
such that for any triple $(P_0,u_0,v_0) \in W$ and $s\in (-a,a)$, 
the triples $\T_s^i(P_0,u_0,v_0)$ are defined for $i= -2,-1,\ldots,3,4$
and depend smoothly on the tuple $(s, P_0,u_0,v_0)$, giving
smooth morphisms $(-a,a)\times W \to TM\times_MTM$.
Given any open neighbourhood $U$ of
$x\in M$, if we choose $a>0$ and $W$ to be sufficiently
small then the above quadrilaterals lie entirely in $U$. 

In terms of the notation introduced above, 
we now restate the Theorem \ref{geodesic gap} in the following
somewhat stronger form, which allows $(P_0,u_0,v_0)$ to vary.

\theorem\label{geodesic gap again}
{\it Let $M$ be a smooth manifold with an affine connection $\nabla$.
Then for any
$C^{\infty}$ function $f$ defined on a neighbourhood $U$ of a point
$x\in M$, and for any $u_x,v_x\in T_xM$, 
the following relations hold in the ring $C^{\infty}((-a,a)\times W)$
for sufficiently small $a>0$ and open neighbourhood $W$ of $(x,u_x,v_x)$
in $TM\times_MTM$, where $(s^n) \subset C^{\infty}((-a,a)\times W)$ 
denotes the principal ideal generated by $s^n$. 

(1) $f(P_2(s)) - f(Q_2(s)) = - s^2T(u_0,{v_0})(f) \mbox{ mod }(s^3)$, 
where $T$ is the torsion tensor of $\nabla$. 

(2) $f(P_4(s)) - f(P_0) = - s^2T(u_0,v_0)(f) \mbox{ mod }(s^3)$.


If moreover $\nabla$ is symmetric, then we have

(3) $f(Q_2(s)) - f(P_2(s)) = {1\over 2}s^3(R(u_0,v_0)(u_0+v_0))(f) 
\mbox{ mod }(s^4)$, where $R$ is the curvature tensor of $\nabla$.

(4) $f(P_4(s)) - f(P_0) = {1\over 2}s^3(R(u_0,v_0)(u_0+v_0))(f) 
\mbox{ mod }(s^4)$.

}

\rm

{\bf Strategy of the proof.} It is enough to take the function $f$ above 
to be a coordinate function $x^i$ w.r.t a smooth local chart 
$(U, x^1,\ldots,x^d)$ around $x$.
We will fix $(P_0,u_0,v_0)$ to begin with, and 
make a laborious but straight-forward computation of the 
Taylor series in $s$ for the function 
$x^i(Q_2(s)) - x^i(P_2(s))$, which will yield the statements (1) and (3). 
By Lemma \ref{Taylor}(2), these imply the full forms of (1) and (3)
(with variable $(P_0,u_0,v_0)$).
Note that we have already proved (1) as a consequence of 
Lemma \ref{meaning of bracket} and Theorem \ref{frame flow bracket}, 
so this gives another proof of the same.

Finally, we employ a certain trick involving 
an identity between two kinds of gaps and a Taylor series expansion
(which is a somewhat more complicated
version of the trick that we use in the proof of 
Lemma \ref{meaning of bracket again}), 
to show that the full form of (1) (with variable $(P_0,u_0,v_0)$) implies 
the form of (2) where $(P_0,u_0,v_0)$ is fixed,
and the full form of (3) implies the form of (4) 
where $(P_0,u_0,v_0)$ is fixed. Again, the full form of (2) and (4)
follows by Lemma \ref{Taylor}(2).
The reverse implications also follow from similar arguments.
The Theorem \ref{geodesic gap} stated
in the Introduction is an immediate consequence of the 
Theorem \ref{geodesic gap again}.


{\bf Proofs of \ref{geodesic gap again}.(1) and \ref{geodesic gap again}.(3).}
To begin with, we fix the triple $(P_0,u_0,v_0)$. 
Let $\gamma_1$ be the geodesic on $M$ with affine parameter $s$,
normalized by $\gamma_1(0) = P_0$ and ${d\gamma_1\over ds}(0) = u_0$.
Let $u_1(s), v_1(s) \in T_{\gamma_1(s)}M$ be the parallel transports of 
${u_0}, {v_0}\in T_{P_0}M$ along $\gamma_1$. 
Let $\gamma_{2,s}$ be the geodesic 
on $M$ with affine parameter $t$,
normalized by $\gamma_{2,s}(0) = \gamma_1(s)$ 
and ${d\gamma_{2,s}\over dt}(0) = v_1(s)$. 
For $s,t$ in a small domain $D = (-a,a)\times (-a,a)$ in $\R^2$, let
$$\P(s,t) = \gamma_{2,s}(t) \in M$$
which defines a smooth function $\P : D \to M$. In particular, 
$\P(s,0) = \gamma_1(s)$.

Let $U$ be a coordinate chart around $P_0$ with coordinates
$(x^1,\ldots, x^d)$, with $P_0 = (x_0^1,\ldots, x_0^d)$.
We choose $a>0$ to be small enough such that the image of 
$\P : D \to M$ lies in $U$. 
Let $x^i(s,t)$ denote the coordinates of $\P(s,t)$.
We use the notation 
$\Gamma^i_{jk}(0) = \Gamma^i_{jk}(P_0)$
and $\Gamma^i_{jk}(s,t) = \Gamma^i_{jk}(\P(s,t))$.
The point $\gamma_1(s)$ has coordinates 
$$x_1^i(s) = x^i(s,0).$$
In terms of the summation convention, we have
${u_0} = u_0^i{\pa \over \pa x^i}$, 
${v_0} = v_0^i{\pa \over \pa x^i}$, 
$u_1(s) = u_1^i{\pa \over \pa x^i}$ and 
$v_1(s) = v_1^i{\pa \over \pa x^i}$,
where $u_1^i$ and $v_1^i$ are the functions of $s$ that are 
the coefficients of $u_1(s)$ and $v_1(s)$.

The curve $\gamma_1$ satisfies the geodesic equation, hence
$${d^2x_1^i\over ds} = {du_1^i\over ds} = - \Gamma^i_{jk}(s)u_1^j(s)u_1^k(s)$$
along $\gamma_1$. 
Also, as $v_1(s)$ is parallel transported along $\gamma_1$, we have 
$${dv_1^i\over ds} = - \Gamma^i_{jk}(s)v_1^j(s)u_1^k(s).$$

The steps in the following calculation 
are obtained by taking Taylor expansions in $s$ 
and making various substitutions. 
\begin{align*}
x_1^i(s) =\ & 
x^i_0 + s {dx_1^i \over ds}(0) + {s^2\over 2}{d^2x_1^i \over ds^2}(0) 
+ {s^3 \over 3!}{d^3x_1^i \over ds^3}(0) \mbox{ mod }(s^4). \\
{dx_1^i \over ds} = \ & u_1^i \mbox{ as }\gamma \mbox{ is a geodesic
                     with tangent }u_1.\\
{d^2x_1^i \over ds^2} = {du_1^i \over ds} =\ & - \Gamma^i_{jk}u_1^ju_1^k
 \mbox{ by the geodesic equation}.\\
{d^3x_1^i \over ds^3} = \ & 
-\Gamma^i_{jk,{\ell}}u_1^ju_1^ku_1^{\ell} + \Gamma^i_{jk}\Gamma^j_{ab}u_1^au_1^bu_1^k 
+ \Gamma^i_{jk}\Gamma^k_{ab}u_1^au_1^bu_1^j. \\
\end{align*}
It follows that in the ring $C^{\infty}(-a,a)$ we have the relation
\begin{align*}
x_1^i(s) =\ & x^i_0 + s u_0^i - {s^2\over 2}\Gamma^i_{jk}(0)u_0^ju_0^k \\
         & + {s^3 \over 3!}(-\Gamma^i_{jk,{\ell}}(0)u_0^ju_0^ku_0^{\ell} 
           + \Gamma^i_{jk}(0)\Gamma^j_{ab}(0)u_0^au_0^bu_0^k 
           + \Gamma^i_{jk}(0)\Gamma^k_{ab}(0)u_0^au_0^bu_0^j)  \mbox{ mod }(s^4). \\
\end{align*}
By Taylor expansion and substitutions, we have the following. 
\begin{align*}
v_1^i(s) =\ &  v_0^i + s {dv_1^i\over ds}(0) + 
{s^2\over 2} {d^2v_1^i\over ds^2}(0)
  \mbox{ mod }(s^3).\\
{dv_1^i\over ds} =\ &  - \Gamma^i_{jk}v_1^ju_1^k
\mbox{ as }v_1(s) \mbox{ is parallel transported along }\gamma_1.\\
{d^2v_1^i\over ds^2} = \ & 
-\Gamma^i_{jk,l}v_1^ju_1^ku_1^l + \Gamma^i_{jk}\Gamma^j_{mn}v_1^mu_1^nu_1^k 
+ \Gamma^i_{jk}\Gamma^k_{mn}v_1^ju_1^mu_1^n.\\
\end{align*}
Hence we get
\begin{align*}
v_1^i(s) = \ & v_0^i + s (- \Gamma^i_{jk}(0)v_0^ju_0^k)\\ 
 & + {s^2\over 2}( -\Gamma^i_{jk,l}(0)v_0^ju_0^ku_0^l 
   + \Gamma^i_{jk}(0)\Gamma^j_{mn}(0)v_0^mu_0^nu_0^k 
   + \Gamma^i_{jk}(0)\Gamma^k_{mn}(0)v_0^ju_0^mu_0^n)\\
 &  \mbox{ mod }(s^3) \mbox{ in the ring }C^{\infty}(-a,a).\\
\end{align*}
Next, we consider 
the geodesic $\gamma_{2,s}(t)$, which we write as
$\gamma_2(t)$ or just $\gamma_2$ when $s$ is suppressed from the notation.
We will denote the coordinates of $\gamma_2(t)$ by $x_2^i(t)$, leaving out
the mention of $s$. 
The curve $\gamma_{2,s}(t)$ is defined by $\gamma_{2,s}(0) = P_1(s)$, and
$${\pa \gamma_{2,s}^i\over \pa t}(s,0) = v_1^i(s)$$
which we will simply write as 
${dx_2^i\over dt}(0) = v_1^i$ by suppressing $s$ from the notation.

The steps in the following calculation 
are obtained by Taylor expansion in 
$t$ and making substitutions. The coefficients of $t^n$ 
terms are functions of $s$, belonging to $C^{\infty}(-a,a)$. 
The various relations are 
in the ring $C^{\infty}(D)$ 
modulo certain ideals $J\subset C^{\infty}(D)$.
We suppress the mention of $s$ and $t$ in some 
places for brevity of notation, when these are understood. 
\begin{align*}
x_2^i(t) = \ &  x^i_1 + t {dx_2^i \over dt}(0) + 
{t^2\over 2}{d^2x_2^i \over dt^2}(0) + {t^3\over 3!}{d^3x_2^i \over dt^3}(0) 
 \mbox{ mod }t^4.\\
{dx_2^i \over dt}(0)  = \ &  v_1^i(s)\\
   = \ & v_0^i + s (- \Gamma^i_{jk}(0)v_0^ju_0^k)\\ 
       & + {s^2\over 2}( -\Gamma^i_{jk,l}(0)v_0^ju_0^ku_0^l 
          + \Gamma^i_{jk}(0)\Gamma^j_{mn}(0)v_0^mu_0^nu_0^k 
            + \Gamma^i_{jk}(0)\Gamma^k_{mn}(0)v_0^ju_0^mu_0^n)\\
 &  \mbox{ mod }s^3.\\
{d^2x_2^i \over dt^2} = \ & -\Gamma^i_{jk} {dx_2^j \over dt}{dx_2^k \over dt}
\mbox{ by the geodesic equation}.\\
{d^2x_2^i \over dt^2}(0) = \ & -\Gamma^i_{jk}(s,0)v_1^j(s)v_1^k(s). \\
\Gamma^i_{jk}(s,0) =  \ & \Gamma^i_{jk}(0) + s\Gamma^i_{jk,\ell}(0)u_0^{\ell} 
 \mbox{ mod }s^2.\\
{d^2x_2^i \over dt^2}(0) =  \ & 
           -(\Gamma^i_{jk}(0) + s\Gamma^i_{jk,\ell}(0)u_0^{\ell})
            (v_0^j - s\Gamma^j_{ab}(0)v_0^au_0^b)
            (v_0^k - s\Gamma^k_{ab}(0)v_0^a u_0^b)
 \mbox{ mod }s^2,\\
=  \ & -\Gamma^i_{jk}(0)v_0^jv_0^k  \\
     &   -s(\Gamma^i_{jk,\ell}(0)u_0^{\ell}v_0^jv_0^k
          -\Gamma^i_{jk}(0) \Gamma^j_{ab}(0)v_0^au_0^bv_0^k
           -\Gamma^i_{jk}(0)\Gamma^k_{ab}(0)v_0^a u_0^bv_0^j)
              \mbox{ mod }s^2.\\
{d^3x_2^i \over dt^3} = \ & 
- \Gamma^i_{jk,\ell}{dx_2^{\ell}\over dt}{dx_2^j \over dt}{dx_2^k \over dt}
+ \Gamma^i_{jk}\Gamma^j_{ab} {dx_2^a \over dt}{dx_2^b \over dt}{dx_2^k \over dt}
+ \Gamma^i_{jk}\Gamma^k_{ab} {dx_2^a \over dt}{dx_2^b \over dt}{dx_2^j \over dt}.
\end{align*}

We now evaluate the last equation at $t=0$, that is, 
at $(s,0)\in D$, modulo the principal ideal $(s)$. 
Note that ${dx_2^j \over dt}(s,0) = v_1^j(s)$, so 
${dx_2^j \over dt}(s,0) = v_0^j \mbox{ mod }s$. 
Similarly, $\Gamma^i_{jk}(s,0) = \Gamma^i_{jk}(0) \mbox{ mod }s$,
and $\Gamma^i_{jk,\ell}(s,0) = \Gamma^i_{jk,\ell}(0) \mbox{ mod }s$.
Hence, 
\begin{align*}
{d^3x_2^i \over dt^3}(0) = \ & 
- \Gamma^i_{jk,\ell}(0)v_0^{\ell}v_0^jv_0^k
+ \Gamma^i_{jk}(0)\Gamma^j_{ab}(0) 
v_0^av_0^bv_0^k
+ \Gamma^i_{jk}(0)\Gamma^k_{ab}(0)
v_0^av_0^bv_0^j\mbox{ mod }s
\end{align*}
Hence by substitutions, we get the following equations in the ring 
$C^{\infty}(D)$ modulo the ideal 
$(s,t)^4=(s^4,s^3t,s^2t^2,st^3,t^4) \subset C^{\infty}(D)$.
\begin{align*}
x^i_2(t) = \ &  x^i_1 + t {dx_2^i \over dt}(0) + 
{t^2\over 2}{d^2x_2^i \over dt^2}(0) + {t^3\over 3!}{d^3x_2^i \over dt^3}(0) 
 \mbox{ mod }(t^4)\\
= \ &  x^i_0 + s u_0^i + t v_0^i - {s^2\over 2}\Gamma^i_{jk}(0)u_0^ju_0^k 
- {t^2\over 2} (\Gamma^i_{jk}(0)v_0^jv_0^k) + ts (- \Gamma^i_{jk}(0)v_0^ju_0^k)\\ 
         & + {s^3 \over 3!}(-\Gamma^i_{jk,{\ell}}(0)u_0^ju_0^ku_0^{\ell} 
           + \Gamma^i_{jk}(0)\Gamma^j_{ab}(0)u_0^au_0^bu_0^k 
           + \Gamma^i_{jk}(0)\Gamma^k_{ab}(0)u_0^au_0^bu_0^j)\\
         & + {s^2t\over 2}( -\Gamma^i_{jk,l}(0)v_0^ju_0^ku_0^l 
                 + \Gamma^i_{jk}(0)\Gamma^j_{mn}(0)v_0^mu_0^nu_0^k 
                 + \Gamma^i_{jk}(0)\Gamma^k_{mn}(0)v_0^ju_0^mu_0^n) \\
         &       - {st^2\over 2}((\Gamma^i_{jk,\ell}(0)u_0^{\ell}v_0^jv_0^k
                      -  \Gamma^i_{jk}(0) \Gamma^j_{ab}(0)v_0^au_0^bv_0^k
                      -  \Gamma^i_{jk}(0)\Gamma^k_{ab}(0)v_0^a u_0^bv_0^j)\\
         & + {t^3\over 3!} (- \Gamma^i_{jk,\ell}(0)v_0^{\ell}v_0^jv_0^k
           + \Gamma^i_{jk}(0)\Gamma^j_{ab}(0)v_0^av_0^bv_0^k
           + \Gamma^i_{jk}(0)\Gamma^k_{ab}(0)v_0^av_0^bv_0^j) \\
         & \mbox{ modulo } (s,t)^4.
\end{align*}

The point $P_2(s)$ is the point $\P(s,s)$. Hence 
putting $t=s$ in the above expression, 
renaming dummy indices and collecting terms, we see that 
the coordinates of $P_2(s)$ as a function of $s\in (-a,a)$ 
are given in the ring $C^{\infty}(-a,a)$ modulo the ideal $(s^4)$ by
\begin{align*}
P_2(s)^i = \ 
 & x^i_0  
   + s(u_0^i +v_0^i) 
   - {s^2\over 2}\Gamma^i_{jk}(0)(u_0^ju_0^k + v_0^jv_0^k) 
   - s^2\Gamma^i_{jk}(0)v_0^ju_0^k \\
 & + {s^3 \over 3!}
(- \Gamma^i_{pq,r}(0) + \Gamma^i_{jr}(0)\Gamma^j_{pq}(0) 
+ \Gamma^i_{rk}(0)\Gamma^k_{pq}(0))(u_0^pu_0^qu_0^r+ v_0^pv_0^qv_0^r)  \\
& + {s^3\over 2}
         (-\Gamma^i_{bc,a}(0)+\Gamma^i_{jb}(0) \Gamma^j_{ca}(0)   
           + \Gamma^i_{ck}(0)\Gamma^k_{ba}(0)) 
(u_0^a u_0^bv_0^c + u_0^a v_0^bv_0^c)\\
         &     \mbox{ mod }(s^4).
\end{align*}
A similar calculation, in which the roles of $u_0$ and $v_0$ are interchanged,
gives the coordinates $Q_2(s)^i$ of the point $Q_2(s)$.
Subtracting the two expressions, 
we get in the ring $C^{\infty}(-a,a)$ 
the following relations modulo the ideal $(s^3)$ 
\begin{align*}
P_2(s)^i - Q_2(s)^i = \ & 
s^2(-\Gamma^i_{jk}(0)v_0^ju_0^k + \Gamma^i_{jk}(0)u_0^jv_0^k )\\ 
= \ & s^2(\Gamma^i_{pq}(0) -\Gamma^i_{qp}(0) )u_0^pv_0^q\\
= \ & -s^2T_{pq}(0)u_0^pv_0^q \mbox{ where }T_{pq}\mbox{ is the torsion tensor}\\
= \ & -s^2T(u_0,v_0)^i ~~\mbox{ mod }(s^3).
\end{align*}
This proves the part (1) of 
Theorem \ref{geodesic gap again}.
Now suppose that $T(u_0,v_0) = 0$, or in coordinate terms,
$(-\Gamma^i_{jk} + \Gamma^i_{kj})u_0^jv_0^k = 0$ for all $i$.
Then in the calculation of $P_2(s)^i - Q_2(s)^i $ modulo $(s^4)$
a large number of terms cancel and we get
$$P_2(s)^i - Q_2(s)^i = -{s^3\over 2}
(\Gamma^i_{pr,q}-\Gamma^i_{pq,r} +\Gamma^i_{qj}\Gamma^j_{pr}
- \Gamma^i_{rj} \Gamma^j_{pq})
u_0^qv_0^r(u_0^p+ v_0^p)
\mbox{ mod } (s^4)$$
Note that  
$$\Gamma^i_{pr,q}-\Gamma^i_{pq,r} +\Gamma^i_{qj}\Gamma^j_{pr}
- \Gamma^i_{rj} \Gamma^j_{pq} = R^i_{pqr}$$
is just the coordinate form of the Riemann curvature tensor.
Hence the above reads
$$P_2(s)^i - Q_2(s)^i = -{s^3\over 2}R^i_{pqr}(0)
u_0^qv_0^r(u_0^p+ v_0^p)
\mbox{ mod } (s^4).$$
Hence we have proved 
the statement (3) of the Theorem \ref{geodesic gap again}.


{\bf Proof that \ref{geodesic gap again}.(1) $\Rightarrow$ 
\ref{geodesic gap again}.(2) and \ref{geodesic gap again}.(3) $\Rightarrow$ 
\ref{geodesic gap again}.(4).} 
As the first step, we define two kinds 
of geodesic gap functions $G_I$ and $G_{II}$. With hypothesis and notation as
in the statement of Theorem \ref{geodesic gap again}, we put
$$G_I(P_0,u_0,v_0,f,s) = f(P_4(s)) - f(P_0) \mbox{ and }
G_{II}(P_0,u_0,v_0,f,s) = f(P_2(s)) - f(Q_2(s)).$$
For a given $f$, we regard both $G_I$ and $G_{II}$ as smooth functions
$(-a,a) \times W \to \R$. 


By a scaling argument, in order to prove the result for a triple $(P,u,v)$,
it is enough to verify it for $(P, \lambda u, \lambda v)$ for 
all $\lambda$ in a neighbourhood of $0$. Hence it is enough to take 
the point $(x,u_x,v_x) \in TM\times_MTM$ to be of the form $(x,0,0)$.
Then note that $\T_0(x,0,0) = \T_0^{-1}(x,0,0) = (x,0,0)$. 
Hence by continuity, $(x,0,0)$ has an open neighbourhood $V$, and there
is some $0< b < a$ such that 
$\T_s^i(P_0,u_0,v_0) \in W$ for all $-2\le i \le 4$, whenever $s\in (-b,b)$
and $(P_0,u_0,v_0) \in V$. 
For a given $(P_0,u_0,v_0)$, recall that 
$$\T_s^2(P_0,u_0,v_0) = (P_2(s), -u_2(s), -v_2(s))$$
where $u_2(s)$ and $v_2(s)$ are respectively the 
parallel transports of $u_0$ and $v_0$ along the leg $P_0P_1(s)$
followed by the leg $P_1(s)P_2(s)$ of the geodesic quadrilaterals.
Now consider the map 
$(-b,b) \to W$ that sends 
$$\tau \mapsto (P_{0,\tau}, u_{0,\tau}, v_{0,\tau})
= (P_2(\tau), - u_2(\tau), -v_2(\tau))$$
Composing the gap function $G_{II}$ 
with this function $(-b,b) \to W$, we get a {\bf parametric
form} of $G_{II}$, defined by
$$(s,\tau) \mapsto G_{II}(P_0,u_0,v_0,f,s,\tau) 
= G_{II}(P_{0,\tau},u_{0,\tau},v_{0,\tau},f,s)$$
which for a given $f$ defines a smooth function
$(-a,a) \times (-b,b) \to \R$.
Hence applying the homomorphism $C^{\infty}((-a,a)\times W) \to 
C^{\infty}((-a,a) \times (-b,b))$ to the relation 
of Theorem \ref{geodesic gap again}.(1) and (3), we get the following
relations in $C^{\infty}((-a,a)\times (-b,b))$.
\begin{align*}
G_{II}(P_{0,\tau},u_{0,\tau},v_{0,\tau},f,s) 
= \ & -s^2T(u_{0,\tau},v_{0,\tau})(f) 
\mbox{ mod }s^3C^{\infty}((-a,a)\times (-b,b)),\\
\mbox{and when }T \equiv 0, ~~~~~~~~~~~~~~~~~~~~ & \\ 
G_{II}(P_{0,\tau},u_{0,\tau},v_{0,\tau},f,s) 
= \ & {s^3\over 2}R(u_{0,\tau},v_{0,\tau})(u_{0,\tau}+v_{0,\tau})(f) \\
& ~~~~~~~~~ \mbox{ mod }s^4C^{\infty}((-a,a)\times (-b,b)).
\end{align*}

Note that at $s=0$, 
the smooth function $T(-u_2(s),- v_2(s))(f)$
takes the value $T(u_0,v_0)(f)$, so we have
$$T(-u_2(s),- v_2(s))(f) =  T(u_0,v_0)(f) \mbox{ mod }(s).$$
Similarly,   
$$R(-u_2(s), -v_2(s))(-u_2(s) -v_2(s))(f)= - R(u_0,v_0)(u_0+v_0)(f) 
\mbox{ mod }(s).$$

The diagonal $\Delta \subset (-a,a)\times (-b,b)$ is defined by
the equation $s = \tau$. This gives a closed imbedding 
$$(-b,b)\hra (-a,a) \times  (-b,b) :
\tau \mapsto (\tau,\tau).$$ 
Restricting functions on 
$(-a,a)\times  (-b,b)$ to $(-b,b)$ 
under the above diagonal imbedding  defines a homomorphism
$$C^{\infty}((-a,a)\times (-b,b))\to C^{\infty}(-b,b)$$
under which $\tau \mapsto s$ and $s\mapsto s$. 
Under this homomorphism, the extension of the principal ideal
$sC^{\infty}((-a,a)\times (-b,b))$ is the principal ideal
$sC^{\infty}(-b,b)$. Hence putting $\tau = s$ in the 
relation 
$$G_{II}(P_{0,\tau},u_{0,\tau},v_{0,\tau},f,s) 
= -s^2T(u_{0,\tau},v_{0,\tau})(f) \mbox{ mod }s^3C^{\infty}((-a,a)\times (-b,b))$$
we get 
\begin{eqnarray*}
G_{II}(P_2(s),-u_2(s),v_2(s),f,s) 
& = & 
-s^2T(-u_2(s) ,-v_2(s))(f) \mbox{ mod }s^3C^{\infty}(-b,b)\\
& = & -s^2T(-u_0 ,-v_0)(f) \mbox{ mod }s^3C^{\infty}(-b,b)\\
& = & -s^2T(u_0 ,v_0)(f) \mbox{ mod }s^3C^{\infty}(-b,b) ~~~\ldots {\bf (*)}
\end{eqnarray*}

Now suppose that $T\equiv 0$, so that we have
$$G_{II}(P_{0,\tau},u_{0,\tau},v_{0,\tau},f,s) = 
{s^3\over 2}R(u_{0,\tau},v_{0,\tau})(u_{0,\tau}+v_{0,\tau})(f) 
\mbox{ mod }s^4C^{\infty}((-a,a)\times (-b,b))$$
Restricting under the diagonal imbedding 
$(-b,b)\hra (-a,a)\times (-b,b)$, we get
the following relations in the ring $C^{\infty}(-b,b)$: 
\begin{eqnarray*}
G_{II}(P_2(s),-u_2(s),-v_2(s),f,s) 
& = & 
{s^3\over 2}R(-u_2(s),-v_2(s))(-u_2(s) - v_2(s))(f) \mbox{ mod }(s^4),\\
& = & {s^3\over 2}R(-u_0, -v_0)(-u_0 - v_0)(f)
 \mbox{ mod }(s^4),\\
& = & -{s^3\over 2}R(u_0, v_0)(u_0+ v_0)(f)
 \mbox{ mod }(s^4)~~~\ldots {\bf (**)}
\end{eqnarray*}

It now remains to deduce the statements about $f(P_4(s)) - f(P_0)$.
For this, recall that we have defined the other gap function $G_I$ by  
$G_I(P_0,u_0,v_0,f,s) = f(P_4(s)) - f(P_0)$.
As $G_{II}$ was defined by
$G_{II}(P_0,u_0,v_0,f,s) = f(P_2(s)) - f(Q_2(s))$,
we have the obvious identity 
$$G_I(P_0,u_0,v_0,f,s) =  G_{II}(P_2(s),-u_2(s),-v_2(s),f,s).$$
We have already proved (see equation $(*)$ above)
that the right hand side is congruent
modulo $(s^3)$ to $- s^2T(u_0,v_0)$. Therefore we have 
$$G_I(P_0,u_0,v_0,f,s) = - s^2T(u_0,v_0)\mbox{ mod }(s^3)$$
in $C^{\infty}(-b,b)$ which completes the proof of 
Theorem \ref{geodesic gap again}.(2) for a fixed $(P_0,u_0,v_0)$, 
and hence also for the general case 
of a variable $(P_0,u_0,v_0)$ 
by applying Lemma \ref{Taylor} 
as already explained. 

By equation $(**)$, 
if $T\equiv 0$ then 
$G_{II}(P_2(s),-u_2(s),-v_2(s),f,s)$ is congruent 
modulo $(s^4)$ to $-{s^3\over 2}R(u_0, v_0)(u_0+ v_0)(f)$.
Therefore we have 
$$G_I(P_0,u_0,v_0,f,s) = {s^3\over 2}R(u_0, v_0)(u_0+ v_0)(f)
\mbox{ mod }(s^4)$$
in $C^{\infty}(-b,b)$.
This completes the proof of Theorem \ref{geodesic gap again}.(4)
for a fixed $(P_0,u_0,v_0)$, 
and hence also for the general case 
of a variable $(P_0,u_0,v_0)$  
as already explained.

\bigskip

{\bf Acknowledgement} I thank Professor Rajaram Nityananda 
for his interesting question that led to this work. 
I thank Dr Ananya Chaturvedi for her help with the derivative calculations,
producing the figures, and for a careful proof reading. 


\bigskip

\centerline{\bf References}

[1] K. Nomizu: {\it `Lie Groups and Differential Geometry'}, 
Publ. Math. Soc. Japan, 1956.

\bigskip

\bigskip

School of Mathematics, Tata Institute of Fundamental Research,
Homi Bhabha Road, Mumbai 400 005, India. $~~~$ email: {\tt nitsure@gmail.com}

\bigskip

\bigskip


\end{document}